\documentclass[12pt]{article}

\usepackage{amsmath,amssymb,amsthm}

\iftrue \makeatletter
\@ifundefined{every@math@size}{}{%
\addto@hook\every@math@size{\dch@scr@hook} 
\def\dch@scr@adjust{\@ifundefined{dch@sizet\f@size}%
  {\expandafter\dch@set@script\csname dch@sizet\f@size\endcsname}%
  {\csname dch@sizet\f@size\endcsname}} 
\def\dch@set@script#1{\begingroup %
  \frozen@everymath{}%
  \let#1\@empty \let\dch@do@one\relax 
  \dch@set@one\scriptscriptstyle\scriptscriptfont#1\ssf@size 
  \dch@set@one\scriptstyle\scriptfont#1\sf@size 
  \dch@set@one\textstyle\textfont#1\f@size 
  \endgroup #1} %
\def\dch@set@one#1#2#3#4{%
  \@ifundefined{dch@size#4}%
   {\expandafter\xdef\csname dch@size#4\endcsname{%
      \fontdimen13\the#2\tw@=\the\fontdimen13#2\tw@ 
      \fontdimen14\the#2\tw@=\the\fontdimen14#2\tw@ 
      \fontdimen15\the#2\tw@=\the\fontdimen15#2\tw@ 
      \fontdimen16\the#2\tw@=\the\fontdimen16#2\tw@ 
      \fontdimen17\the#2\tw@=\the\fontdimen17#2\tw@}%
  }{\csname dch@size#4\endcsname}%
  \setbox\z@\hbox{$#1H_2$}\@tempdima\dp\z@ 
  \setbox\z@\hbox{$#1H_2^{+\vrule \@height 1em}$}%
   \ifdim\@tempdima<\dp\z@ 
    \advance\@tempdima\dp\z@ \divide\@tempdima\tw@ %
    \@tempdimb\dp\z@ \advance\@tempdimb-\@tempdima %
    \advance\@tempdimb\ht\z@ \advance\@tempdimb-1em %
    \xdef#3{#3\dch@do@one#2{\the\@tempdimb}{\the\@tempdima}}%
  \fi} 
\def\dch@do@one#1#2#3{%
  \fontdimen13#1\tw@#2\relax 
  \fontdimen14#1\tw@\fontdimen13#1\tw@ \fontdimen15#1\tw@\fontdimen13#1\tw@ 
  \fontdimen\sixt@@n#1\tw@#3\fontdimen17#1\tw@\fontdimen\sixt@@n#1\tw@}%
\let\dch@scr@hook\dch@scr@adjust 
\ifx\glb@currsize\f@size \dch@scr@adjust \fi}%
\makeatother \fi

\date{}

\normalsize

\newtheorem{thm}{Theorem}

\newtheorem{lemma}{Lemma}

\def\dfrac#1#2{\lower0.15ex\hbox{\large$\frac{#1}{#2}$}}
\def\({\bigl(}
\def\){\bigr)}

\def\sumpp{\sum\nolimits'}
\def\sumjkmn{\sum_{j=1}^m\sum_{k=1}^n}
\def\jsum{\sum_{j=1}^{m-1}}
\def\ksum{\sum_{k=1}^{n-1}}
\def\Amin{A_{\mathrm{min}}}

\def\b{{\mathord{\bullet}}}  
\def\c{{\mathord{\ast}}}       

\def\Vd{V_{\rm d}^{}}
\def\Vdinv{V_{\rm d}^{-1}}
\def\Vnd{V_{\rm nd}^{}}

\def\A{{\cal A}}

\def\M{{\cal M}}

\def\Q{{\cal Q}}
\def\R{{\cal R}}
\def\S{{\cal S}}

\def\W{{\cal W}}
\def\X{{\cal X}}
\def\Y{{\cal Y}}

\let\eps=\varepsilon

\let\t=\theta
\let\p=\phi
\def\Deltait{\mathit{\Delta}}
\def\Thetait{\mathit{\Theta}}
\def\Gammait{\mathit{\Gamma}}
\def\varnu{{\dot\nu}}

\def\iter#1{^{(#1)}}

\def\svec{\boldsymbol{s}}
\def\tvec{\boldsymbol{t}}
\def\xvec{\boldsymbol{x}}

\def\zvec{\boldsymbol{z}}

\def\betavec{\boldsymbol{\betavec}}
\def\thetavec{\boldsymbol{\theta}}
\def\phivec{\boldsymbol{\phi}}
\def\sigmavec{\boldsymbol{\sigma}}
\def\tauvec{\boldsymbol{\tau}}
\def\xivec{\boldsymbol{\xi}}
\def\zetavec{\boldsymbol{\zeta}}

\def\avtheta{\bar{\thetavec}}
\def\avphi{\bar{\phivec}}
\def\avx{\bar{\xvec}}

\def\Chat{\hat C}
\def\thetahatvec{\hat{\thetavec}}
\def\phihatvec{\hat{\phivec}}
\def\that{\hat\theta}  
\def\phat{\hat\phi}    
\def\varthetahatvec{\hat{\boldsymbol{\vartheta}}}
\def\varphihatvec{\hat{\boldsymbol{\varphi}}}

\def\mw#1{\hat{#1}}
\def\mwA{\mw{A}}
\def\mwB{\mw{B}}
\def\mwC{\mw{C}}

\def\mwE{\mw{E}}
\def\mwF{\mw{F}}

\def\mwJ{\mw{J}}
\def\mwZ{\mw{Z}}

\def\mwa{\mw{a}}
\def\Mom{\eta}
\def\Momj{\Mom_j}

\def\OO{\widetilde{O}}

\def\Re{\operatorname{Re}}
\def\Im{\operatorname{Im}}

\def\Prob{\operatorname{Prob}}
\def\Reals{{\mathbb{R}}}
\def\expect{\operatorname{\mathbb E}}
\def\abs#1{\lvert#1\rvert} \let\card=\abs
\DeclareMathOperator{\tr}{tr}
\def\nicebreak{\vskip 0pt plus 50pt\penalty-300\vskip 0pt plus -50pt }
\def\ssqrt#1{\sqrt{\vrule width0pt height1.3ex depth0.4ex
  \smash[b]{#1}}} 

\begin{document}

\title {Asymptotic enumeration of dense 0-1 matrices\\ with
        specified line sums}

\author{
E.~Rodney~Canfield\thanks
 {Research supported by the NSA Mathematical Sciences Program.} \\
\small Department of Computer Science\\[-0.8ex]
\small University of Georgia\\[-0.8ex]
\small Athens, GA 30602, USA\\[-0.3ex]
\small\texttt{ercanfie@uga.edu}
\and
Catherine~Greenhill\thanks
{Research supported by the UNSW Faculty Research Grants Scheme.}\\
\small School of Mathematics and Statistics\\[-0.8ex]
\small University of New South Wales\\[-0.8ex]
\small Sydney, Australia 2052\\[-0.3ex]
\small\texttt{csg@unsw.edu.au}
\and
\vrule width0pt height3ex
Brendan~D.~McKay\vrule width0pt height2ex\thanks
 {Research supported by the Australian Research Council.}\\
\small Department of Computer Science\\[-0.8ex]
\small Australian National University\\[-0.8ex]
\small Canberra ACT 0200, Australia\\[-0.3ex]
\small\texttt{bdm@cs.anu.edu.au}
}

\maketitle

\begin{abstract}
Let $\svec=(s_1,s_2,\ldots, s_m)$ and
$\tvec = (t_1,t_2,\ldots, t_n)$ be vectors of non-negative integers
with  $\sum_{i=1}^{m}s_i = \sum_{j=1}^n t_j$.  Let
$B(\svec,\tvec)$ be the number of $m\times n$
matrices over $\{ 0,1\}$ with $j$-th row sum equal to $s_j$ for
$1\leq j\leq m$
and $k$-th column sum equal to $t_k$ for $1\leq k\leq n$.  Equivalently,
$B(\svec,\tvec)$ is the number of bipartite graphs with
$m$ vertices  in one part with degrees given by~$\svec$,
and $n$ vertices in the other part with degrees given by~$\tvec$.
Most research on the asymptotics of $B(\svec,\tvec)$ has focused on
the sparse case, where the best result is that of
Greenhill, McKay and Wang (2006).
In~the case of dense matrices, the only precise result is for
the case of equal row sums and equal column sums (Canfield
and McKay, 2005).
This paper extends the analytic methods used by the latter
paper to the case where the row and column sums can vary
within certain limits.  Interestingly, the result can be expressed by the
same formula which holds in the sparse case.
\end{abstract}

\nicebreak
\section{Introduction}\label{s:intro}

Let $\svec=(s_1,s_2,\ldots, s_m)$ and
$\tvec = (t_1,t_2,\ldots, t_n)$ be vectors of positive integers
with  $\sum_{i=1}^{m}s_i = \sum_{j=1}^n t_j$.  Let
$B(\svec,\tvec)$ be the number of $m\times n$
matrices over $\{ 0,1\}$ with $j$-th row sum equal to $s_j$ for $1\leq j\leq m$
and $k$-th column sum equal to $t_k$ for $1\leq k\leq n$.  Equivalently,
$B(\svec,\tvec)$ is the number of labelled bipartite graphs with
$m$ vertices  in one part of the bipartition with degrees given
by $\svec$, and $n$ vertices in the other part of the bipartition
with degrees given by $\tvec$.
Let $s$ be the average value of $s_1,s_2,\ldots,s_m$
and let~$t$ be the average value of $t_1,t_2,\ldots,t_n$.
Define the density $\lambda = s/n = t/m$,
which is the fraction of entries in the matrix which equal 1.

The asymptotic value of $B(\svec,\tvec)$ has been much studied,
especially since the celebrated Gale-Ryser Theorem~\cite{Ryser} that
characterises $(\svec,\tvec)$ such that $B(\svec,\tvec)>0$.
Various authors have considered the \emph{semiregular} case,
where $s_j=s$ for $1\leq j\leq m$ and $t_k=t$ for $1\leq k\leq n$.
Write $B(m,s; n,t)$ for $B(\svec,\tvec)$ in this case.
For the sparse (low-$\lambda$) semiregular case,
the best result is by McKay and Wang~\cite{MW}
who gave an asymptotic expression for $B(m,s; n,t)$ which holds
when $st=o\((mn)^{1/2}\)$. 
In the dense ($\lambda$ not close to 0 or 1) semiregular case,
Canfield and McKay~\cite{CM}
used analytic methods to obtain an asymptotic expression for
$B(m,s;n,t)$ in two ranges: in the first, the matrix is relatively square
and the density is not too close to 0 or 1, while in the second, the matrix
is much wider than high (or vice-versa) but the density is arbitrary.
For the sparse irregular case, the best result is that of Greenhill,
McKay and Wang~\cite{GMW}, who gave an asymptotic expression
for  $B(\svec,\tvec)$ which holds when
$\max\{s_j\}\max\{t_k\} = o\((\lambda mn)^{2/3}\)$.

See~\cite{CM}, \cite{GMW} and~\cite{MW} for a more extensive
historical survey.

The contribution of this paper is to adapt the approach of~\cite{CM} to the
dense irregular case when the matrix is relatively square and the density is
not too close to 0 or 1.
See McKay and Wormald~\cite{MWreg} for the corresponding calculation
for symmetric matrices.

In keeping with these earlier papers,
the asymptotic value of $B(\svec,\tvec)$ can be expressed
by a very nice formula involving binomial coefficients.  We now state our
theorem.

\begin{thm}
\label{bigtheorem}
Let $\svec=\svec(m,n)=(s_1,s_2,\ldots,s_m)$ and
$\tvec=\tvec(m,n)=(t_1,t_2,\ldots,t_n)$ be vectors of positive
integers such that
$\sum_{j=1}^m s_j=\sum_{k=1}^n t_k$ for all~$m,n$.
Define $s=m^{-1}\sum_{j=1}^m s_j$, $t=n^{-1}\sum_{k=1}^n t_k$,
$\lambda=s/n=t/m$ and $A = \tfrac12 \lambda(1-\lambda)$.
For some $\eps>0$, suppose
that $\abs{s_j-s}=O(n^{1/2+\eps})$ uniformly for $1\le j\le m$,
and $\abs{t_k-t}=O(m^{1/2+\eps})$ uniformly for $1\le k\le n$.
Define $R=\sum_{j=1}^m (s_j-s)^2$ and $C=\sum_{k=1}^n (t_k-t)^2$.
Let $a,b>0$ be constants such that $a+b<\tfrac{1}{2}$.
Suppose that
$m$, $n\to \infty$ with $n = o(A^2 m^{1+\eps})$, $m=o(A^2 n^{1+\eps})$
 and
\[ \frac{(1-2\lambda)^2}{8A}\biggl(1 + \frac{5m}{6n}+
   \frac{5n}{6m}\biggr) \leq a\log n. \]
Then, provided $\eps>0$ is small enough, we have
\[
 B(\svec,\tvec) = \binom{mn}{\lambda mn}^{\!\!-1}
\prod_{j=1}^m \binom{n}{s_j} \prod_{k=1}^n \binom{m}{t_k}
  \exp\biggl( -\frac12\Bigl(1-\frac{R}{2Amn}\Bigr)
    \Bigl(1-\frac{C}{2Amn}\Bigr) + O(n^{-b})\biggr).
\]
\end{thm}
\begin{proof}
The proof of this theorem is the topic of the paper; here we will summarize
the main phases and draw their conclusions together. 
The basic idea is to identify $B(\svec,\tvec)$ as a
coefficient in a multivariable generating function and to extract that
coefficient using the saddle-point method.
In Section~2, equation~\eqref{start},
we write $B(\svec,\tvec)=P(\svec,\tvec)I(\svec,\tvec)$,
where $P(\svec,\tvec)$ is a rational expression and $I(\svec,\tvec)$
is an integral in $m+n$ complex dimensions.  Both depend on the
location of the saddle point, which is the solution of some nonlinear
equations.  Those equations are solved in Section~\ref{s:radii},
and this leads to the value of $P(\svec,\tvec)$ in~\eqref{rad17}.
In~Section~\ref{s:integral}, the integral $I(\svec,\tvec)$ is estimated
in a small region $\R'$ defined in~\eqref{RSprime}.
The result is given by Theorem~\ref{Jintegral}
together with~\eqref{IvsJ}.  Finally, in Section~\ref{s:boxing}, it is shown that the integral $I(\svec,\tvec)$ restricted to the exterior
of $\R'$ is negligible.  The present theorem thus follows
from~\eqref{start}, \eqref{rad17},
Theorems~\ref{Jintegral}--\ref{boxing} and \eqref{IvsJ}.
\end{proof}

Note that the error term in the above slightly improves the error term for the semiregular case proved in~\cite{CM}.

Thereom~\ref{bigtheorem} has an instructive
interpretation.  Write it as
$B(\svec,\tvec) = N P_1 P_2 E$, where
\begin{align*}
  N &= \binom{mn}{\lambda mn},\quad
  P_1 = N^{-1} \prod_{j=1}^m \binom{n}{s_j},\quad
  P_2 = N^{-1} \prod_{k=1}^n \binom{m}{t_k},\\
   E &= \exp\biggl( -\frac12\Bigl(1-\frac{R}{2Amn}\Bigr)
    \Bigl(1-\frac{C}{2Amn}\Bigr) + O(n^{-b})\biggr).
\end{align*}
Clearly, $N$ is the number of $m\times n$ binary
matrices with $\lambda mn$ ones.  $P_1$ is the 
probability that a matrix randomly chosen from
this class has row sums $\svec$, while $P_2$ is the
probability of the similar event
of having column sums $\tvec$.
If these two events were independent, we would have $E=1$,
so $E$ can be taken as a measure of their non-independence.
For the case when $\svec$ and $\tvec$ are vectors of
constants, that is, $R=C=0$, Ordentlich and Roth~\cite{Ord}
proved that $E\le 1$.

It is proved in \cite{GMW} that the same formula for $B(\svec,\tvec)$
modulo the error term also holds in the sparse case.
Specifically, it holds with a different
vanishing error term whenever
$\max\{s_j\}\max\{t_k\} = o\((\lambda mn)^{2/3}\)$,
$R+C=O\( (\lambda mn)^{4/3}\)$ and 
$RC=O\( (\lambda mn)^{7/3}\)$.
In~\cite{CM}, evidence is presented that the formula is universal
in the semiregular case ($R=C=0$) and it is tempting to conjecture
that the same is true in the irregular case for a wide range of
$R,C$ values.

\medskip

We will use a shorthand notation for summation
over doubly subscripted variables.  If $x_{jk}$ is a variable for
$1\leq j\leq m$ and $1\leq k\leq n$, then
\begin{align*}
  && x_{j\b} &= \sum_{k=1}^{n} x_{jk},&
   x_{\b k} &= \sum_{j=1}^{m} x_{jk}, &
   x_{\b \b } &= \sum_{j=1}^{m}\sum_{k=1}^{n} x_{jk}, &&\\
 && x_{j\c} &= \ksum x_{jk},&
   x_{\c k} &= \jsum x_{jk},&
   x_{\c \c } &= \jsum \ksum x_{jk}, &&
\end{align*}
for $1\leq j\leq m$ and $1\leq k\leq n$.

Throughout the paper, the asymptotic notation $O(f(m,n))$ refers
to the passage of~$m$ and~$n$ to $\infty$.
We also use a modified notation $\OO(f(m,n))$,
which is to be taken as a shorthand for $O\(f(m,n)n^{O(1)\eps}\)$.
In this case it is important that the $O(1)$ factor is uniform
over~$\eps$ provided $\eps$ is small enough; for example we cannot
write $f(m,n)n^{(\eps^{-1})\eps}$ as $\OO(f(m,n))$ even
though $\eps^{-1}=O(1)$ ($\eps$ being defined as a constant).
Under the assumptions of Theorem~\ref{bigtheorem}, we have
$m=\OO(n)$ and $n=\OO(m)$.  We also have that
$8\le A^{-1}\le O(\log n)$, so $A^{-1}=\OO(1)$.
More generally, $A^{c_1}m^{c_2+c_3\eps}n^{c_4+c_5\eps}
=\OO(n^{c_2+c_4})$ if $c_1,c_2,c_3,c_4,c_5$ are constants.

\nicebreak
\section{Expressing the desired quantity as an integral}\label{s:calculations}

In this section we express $B(\svec,\tvec)$ as a contour integral in
$(m+n)$-dimensional complex space, then begin to estimate its value
using the saddle-point method.

Firstly, notice that $B(\svec,\tvec)$ is the coefficient of
$x_1^{s_1}\cdots x_m^{s_m}\,y_1^{t_1}\cdots y_n^{t_n}$ in the function
\[ \prod_{j=1}^m\prod_{k=1}^n \,(1+x_jy_k).\]
By Cauchy's coefficient theorem this equals
\[ B(\svec,\tvec) = \frac{1}{(2\pi i)^{m+n}} \oint \cdots \oint
  \frac{\prod_{j=1}^m\prod_{k=1}^n
              (1+x_jy_k)}{x_1^{s_1+1}\cdots x_m^{s_m+1}
                 y_1^{t_1+1} \cdots y_n^{t_n+1}} \,
              dx_1 \cdots dx_m\,    dy_1 \cdots dy_n,
\]
where each integral is along a simple closed contour
enclosing the origin anticlockwise.
It~will suffice to take each contour to be a circle;
specifically, we will write
\[ x_j = q_j e^{i\theta_j} \quad \mbox{ and } \quad y_k = r_k e^{i\phi_k}\]
for $1\leq j\leq m$ and $1\leq k\leq n$.  Also define
\[ \lambda_{jk} = \frac{q_jr_k}{1+q_jr_k} \]
for $1\leq j\leq m$ and $1\leq k\leq n$.   Then
$1+x_jy_k=(1+q_jr_k)
 \(1+\lambda_{jk}(e^{i(\theta_j+\phi_k)}-1)\)$, so
\begin{align}
\begin{split}\label{start}
 B(\svec,\tvec) &= \frac{\prod_{j=1}^m\prod_{k=1}^n (1 + q_jr_k)}
{(2\pi)^{m+n}\prod_{j=1}^m q_j^{s_j} \prod_{k=1}^n r_k^{t_k}} \\
  &\kern12mm{}\times
  \int_{-\pi}^{\pi}\!\cdots \int_{-\pi}^\pi \frac{\prod_{j=1}^m\prod_{k=1}^n
    (1 + \lambda_{jk}(e^{i(\theta_j+\phi_k)}-1))}
  {\exp(i\sum_{j=1}^m s_j\theta_j + i\sum_{k=1}^n t_k\phi_k)}
   \, d\thetavec d\phivec,
\end{split}
\end{align}
where $\thetavec=(\t_1,\ldots,\t_m)$ and $\phivec=(\p_1,\ldots,\p_n)$.
Write $B(\svec,\tvec) = P(\svec,\tvec) I(\svec,\tvec)$
where $P(\svec,\tvec)$
denotes the factor in front of the integral in~\eqref{start}
and $I(\svec,\tvec)$ denotes the integral.
We will choose the radii
$q_j$, $r_k$ so that there is no linear term in the logarithm
of the integrand of $I(\svec,\tvec)$ when expanded for
small $\thetavec,\phivec$.  This gives the equation
\[
 \sum_{j=1}^m\sum_{k=1}^n\lambda_{jk} (\theta_j+\phi_k)
 - \sum_{j=1}^m s_j\theta_j - \sum_{k=1}^n t_k\phi_k=0.
\]
For this to hold for all $\thetavec, \phivec$, we require
\begin{equation}\label{rad1}
\begin{split}
  \lambda_{j\b} &= s_j \quad (1\le j\le m),\\
  \lambda_{\b k} &= t_k \quad (1\le k\le n).
\end{split}
\end{equation}
In Section~\ref{s:radii} we show that \eqref{rad1} has a solution, and
determine to sufficient accuracy the various functions of the radii,
such as $P(\svec,\tvec)$, that we require.
In Section~\ref{s:integral}
we evaluate the integral $I(\svec,\tvec)$
within a certain region~$\R$ defined in~\eqref{Rdef}.
Section~\ref{s:boxing} contains the
proof that the integral is concentrated within the region $\R$.

\nicebreak
\section{Locating the saddle-point}\label{s:radii}

In this section we solve \eqref{rad1} and derive some of the
consequences of the solution.  As with the whole paper, we
work under the assumptions of Theorem~\ref{bigtheorem}.

Change variables to $\{a_j\}_{j=1}^m$, $\{b_k\}_{k=1}^n$
as follows:
\begin{equation}\label{qrdef}
  q_j = r\frac{1+a_j}{1-r^2a_j},\quad
  r_k = r\frac{1+b_k}{1-r^2b_k},
\end{equation}
where
\[
  r= \sqrt{\frac{\lambda}{1-\lambda}}\;.
\]

Equation \eqref{rad1} is slightly underdetermined, which we will
exploit to impose an additional condition. If $\{q_j\}, \{r_k\}$
satisfy~\eqref{rad1} and $c>0$ is a constant, then $\{cq_j\},
\{r_k/c\}$ also satisfy~\eqref{rad1}.
{}From this we can see that, if there is a solution to~\eqref{rad1}
at all, there is one
for which $\sum_{j=1}^m a_j<0$ and $\sum_{k=1}^n b_k>0$, and
also a solution for which $\sum_{j=1}^m a_j>0$ and
$\sum_{k=1}^n b_k<0$.  It follows from the
Intermediate Value Theorem that there is a solution for which
\begin{equation}\label{rad4}
  n \sum_{j=1}^m a_j = m\sum_{k=1}^n b_k,
\end{equation}
so we will seek a common solution to~\eqref{rad1}
and \eqref{rad4}.

{}From \eqref{qrdef} we find that
\begin{equation}\label{rad2}
  \lambda_{jk}/\lambda = 1 + a_j + b_k + Z_{jk},
\end{equation}
where
\begin{equation}\label{Zjk}
 Z_{jk} = \frac{a_jb_k(1-r^2-r^2a_j - r^2b_k)}
               {1+r^2a_jb_k},
\end{equation}
and that equations~\eqref{rad1} can be rewritten as
\begin{equation}\label{rad3}
\begin{split}
  a_j &= \frac{s_j-s}{\lambda n}- \frac{1}{n}
           \sum_{k=1}^n b_k - \frac{Z_{j\b}}{n}
                      \quad(1\le j\le m),\\
  b_k &= \frac{t_k-t}{\lambda m} - \frac{1}{m}\sum_{j=1}^m a_j
         - \frac{Z_{\b k}}{m} \quad (1\le k\le n).
\end{split}
\end{equation}

\noindent Summing \eqref{rad3} over $j,k$, we find in both cases that
\begin{equation}\label{rad5}
 n\sum_{j=1}^m a_j + m\sum_{k=1}^n b_k = - Z_{\b\b}\,.
\end{equation}
Equations~\eqref{rad4} and~\eqref{rad5} together imply that
\[
     n\sum_{j=1}^m a_j = m\sum_{k=1}^n b_k = -\dfrac12 Z_{\b\b}\,.
\]
Substituting back into~\eqref{rad3}, we obtain
\begin{equation}\label{rad6}
\begin{split}
   a_j &= \mathbb{A}_j(a_1,\ldots,a_m,b_1,\ldots,b_n),\\
   b_k &= \mathbb{B}_k(a_1,\ldots,a_m,b_1,\ldots,b_n),
\end{split}
\end{equation}
for $1\le j\le m$, $1\le k\le n$, where
\begin{equation*}
\begin{split}
  \mathbb{A}_j(a_1,\ldots,a_m,b_1,\ldots,b_n) 
    &= \frac{s_j-s}{\lambda n} - \frac{Z_{j\b}}{n} 
        + \frac{Z_{\b\b}}{2mn}, \\[0.5ex]
  \mathbb{B}_k(a_1,\ldots,a_m,b_1,\ldots,b_n)
   &= \frac{t_k-t}{\lambda m} - \frac{Z_{\b k}}{m} 
        + \frac{Z_{\b\b}}{2mn}.
\end{split}
\end{equation*}

Equation \eqref{rad6} suggests an iteration.
Start with $a\iter0_j=b\iter0_k=0$ for all~$j,k$,  and,
for each $\ell\ge 0$, define
\begin{equation}\label{rad6b}
\begin{split}
   a\iter{\ell+1}_j 
   &= \mathbb{A}_j(a\iter{\ell}_1,\ldots,a\iter{\ell}_m,
     b\iter{\ell}_1,\ldots,b\iter{\ell}_n),\\
   b\iter{\ell+1}_k 
   &= \mathbb{B}_k(a\iter{\ell}_1,\ldots,a\iter{\ell}_m,
   b\iter{\ell}_1,\ldots,b\iter{\ell}_n),
\end{split}
\end{equation}
where $Z_{j\b}=Z\iter{\ell}_{j\b}
=Z_{j\b}(a\iter{\ell}_1,\ldots,a\iter{\ell}_m,
b\iter{\ell}_1,\ldots,b\iter{\ell}_n)$ and similarly for
$Z_{\b k}=Z\iter{\ell}_{\b k}$ and~$Z_{\b\b}=Z\iter{\ell}_{\b\b}$.
We will show that this iteration converges to a solution
of~\eqref{rad6} using a standard contraction-mapping
argument.

Recall that $A^{-1}=O(\log n)$ under the assumptions of
Theorem~\ref{bigtheorem} (which we are adopting throughout).
This implies that $r^2=O(\log n)$.  Therefore,
within the region $\A$ defined by
$\abs{a_j},\abs{b_k}\le n^{-1/3}$ for all~$j,k$,
we have that
\[
\frac{\partial Z_{j,k}}{\partial a_j} = o(m^{-1/4})
\text{\quad and\quad}
\frac{\partial Z_{j,k}}{\partial b_k} = o(n^{-1/4}),
\]
which imply that, in the same region, we have
\begin{align*}
  \frac{\partial \mathbb{A}_j}{\partial a_{j'}} 
    &= \begin{cases} \,o(m^{-1/4}) & (j'=j) \\
                                    \,o(m^{-5/4}) & (j'\ne j),\end{cases} &
  \frac{\partial \mathbb{A}_j}{\partial b_k} &= o(n^{-5/4}), \\
    \frac{\partial \mathbb{B}_k}{\partial a_j} &= o(m^{-5/4}), & 
    \frac{\partial \mathbb{B}_k}{\partial b_{k'}} 
    &= \begin{cases} \,o(n^{-1/4}) & (k'=k) \\
                                    \,o(n^{-5/4}) & (k'\ne k).\end{cases}
\end{align*}
Therefore, by the mean value theorem, we have for $\ell\ge1$ that
\begin{align*}
\begin{split}
  \max_j\, \abs{a\iter{\ell+1}_j -a\iter{\ell}_j}
  + \max_k\,\abs{b\iter{\ell+1}_k -b\iter{\ell}_k}
  &= o(m^{-1/4}) \max_j\, \abs{a\iter{\ell}_j -a\iter{\ell-1}_j} \\
   &\kern5mm{}+o(n^{-1/4}) \max_k\,\abs{b\iter{\ell}_k -b\iter{\ell-1}_k},
\end{split}
\end{align*}
provided $\{a\iter{\ell-1}_j\}\cup\{b\iter{\ell-1}_k\}
\cup\{a\iter{\ell}_j\}\cup\{b\iter{\ell}_k\}\subseteq\A$.

Applying the iteration once, we have
\[
 a\iter1_j=(s_j-s)/(\lambda n)
\text{\quad and\quad}b\iter1_k=(t_k-t)/(\lambda m).
\]
Since $\{a\iter0_j\}$,  $\{b\iter0_k\}$ and
$\{a\iter1_j\}$,  $\{b\iter1_k\}$ lie inside $\frac12\A$,
we find by induction that
$\{a\iter{\ell}_j\}$,  $\{b\iter{\ell}_k\}$ lie in
$\frac{\ell}{\ell+1}\A$  for all $\ell$.
Moreover, the iteration is Cauchy-convergent in the
maximum norm, and the error in stopping at
$\{a\iter{\ell}_j\}, \{b\iter{\ell}_k\}$ is
at most
$\max_j\, \abs{a\iter{\ell}_j -a\iter{\ell-1}_j}
  + \max_k\,\abs{b\iter{\ell}_k -b\iter{\ell-1}_k}$.

When we carry out this iteration, we find that all the
encountered $a\iter{l}_j$ and $b\iter{l}_k$ values are~$\OO(n^{-1/2})$.
It helps to know that the
following approximation holds in that case:
\[
Z_{jk} = (1-r^2)a_jb_k - r^2a_jb_k^2
        -r^2a_j^2b_k - r^2(1-r^2)a_j^2b_k^2 + \OO(n^{-5/2}).
\]

Using the fact that $\sum_{j=1}^m a\iter1_j=0$ and
$\sum_{k=1}^n b\iter1_k=0$, we find that
\begin{align*}
   Z\iter1_{j\b} &= - r^2a\iter1_j\sum_{k=1}^n \(b\iter1_k\)^2 + \OO(n^{-1}), \\
   Z\iter1_{\b k} &= - r^2b\iter1_k\sum_{j=1}^m \(a\iter1_j\)^2 + \OO(n^{-1}), \\
   Z\iter1_{\b\b} &= \OO(1).
\end{align*}
Therefore,
\begin{align*}
  a\iter2_j &= \frac{s_j-s}{\lambda n}
   + \frac{(s_j-s)C}{\lambda^2(1-\lambda)m^2n^2} + \OO(n^{-2}) \quad(1\le j\le m),\\
  b\iter2_k &= \frac{t_k-t}{\lambda m}
   + \frac{(t_k-t)R}{\lambda^2(1-\lambda)m^2n^2} + \OO(n^{-2}) \quad(1\le k\le n).
\end{align*}

Similarly,
\begin{align*}
  Z\iter2_{j\b} &= - r^2a\iter2_j\sum_{k=1}^n \(b\iter2_k\)^2
               - r^2(1-r^2)\(a\iter2_j\)^2\sum_{k=1}^n \(b\iter2_k\)^2 + \OO(n^{-3/2}),\\
  Z\iter2_{\b k} &= - r^2b\iter2_k\sum_{j=1}^m \(a\iter2_j\)^2
                - r^2(1-r^2)\(b\iter2_k\)^2\sum_{j=1}^m \(a\iter2_j\)^2 + \OO(n^{-3/2}),\\
  Z\iter2_{\b\b} &=
         - r^2(1-r^2)\sum_{j=1}^m \(a\iter2_j\)^2\sum_{k=1}^n \(b\iter2_k\)^2 + \OO(n^{-1/2}),
\end{align*}
which gives
\begin{equation}\label{rad8}
\begin{split}
  a\iter3_j &= \frac{s_j-s}{\lambda n}
      + \frac{(s_j-s) C}{\lambda^2(1-\lambda)m^2n^2}
      + \frac{(1-2\lambda)(s_j-s)^2 C}{\lambda^3(1-\lambda)^2m^2n^3}\\
    &\quad{} - \frac{(1-2\lambda) R C}{2\lambda^3(1-\lambda)^2m^3n^3}
          + \OO(n^{-5/2}) \qquad(1\le j\le m),\\
  b\iter3_k &= \frac{t_k-t}{\lambda m}
    + \frac{(t_k-t) R}{\lambda^2(1-\lambda)m^2n^2}
      + \frac{(1-2\lambda)(t_k-t)^2 R}{\lambda^3(1-\lambda)^2m^3n^2}\\
    &\quad{} - \frac{(1-2\lambda) R C}{2\lambda^3(1-\lambda)^2m^3n^3}
          + \OO(n^{-5/2}) \qquad(1\le k\le n).
\end{split}
\end{equation}

Further iterations make no change to this accuracy, so we have that
$a_j = a\iter3_j + \OO(n^{-5/2})$ and $b_k=b\iter3_k + \OO(n^{-5/2})$.
We also have that
\begin{equation}\label{rad9}
\begin{split}
Z_{jk} &=
  \frac{(1-2\lambda)(s_j-s)(t_k-t)}{\lambda^2(1-\lambda)mn}
  - \frac{(s_j-s)(t_k-t)^2}{\lambda^2(1-\lambda)m^2n}
  - \frac{(s_j-s)^2(t_k-t)}{\lambda^2(1-\lambda)mn^2}\\
  &\quad {}
    - \frac{(1-2\lambda)(s_j-s)^2(t_k-t)^2}{\lambda^3(1-\lambda)^2m^2n^2}
  + \frac{(1-2\lambda)(s_j-s)(t_k-t)R}{\lambda^3(1-\lambda)^2m^2n^3}\\
  &\quad {} + \frac{(1-2\lambda)(s_j-s)(t_k-t)C}{\lambda^3(1-\lambda)^2m^3n^2}
         + \OO(n^{-5/2}).
\end{split}
\end{equation}

A sufficient approximation of $\lambda_{jk}$ is given by substituting
\eqref{rad8} and \eqref{rad9} into~\eqref{rad2}.   In~evaluating
the integral $I(\svec,\tvec)$, the following approximations
will be required:
\begin{gather}
\begin{split}\label{rad10}
  \lambda_{jk}(1-\lambda_{jk}) &= \lambda(1-\lambda)
        + \frac{(1-2\lambda)(s_j-s)}{n} + \frac{(1-2\lambda)(t_k-t)}{m}
         - \frac{(s_j-s)^2}{n^2}\\
        &\quad{} - \frac{(t_k-t)^2}{m^2}
        + \frac{(1-6\lambda+6\lambda^2)(s_j-s)(t_k-t)}
               {\lambda(1-\lambda)mn} + \OO(n^{-3/2}),
\end{split} \displaybreak[0]\\[0.5ex]
\begin{split}\label{rad11}
  \lambda_{jk}(1-\lambda_{jk})(1-2\lambda_{jk})
 &= \lambda(1-\lambda)(1-2\lambda)
        + \frac{(1-6\lambda+6\lambda^2)(s_j-s)}{n}\\
        &\qquad{} + \frac{(1-6\lambda+6\lambda^2)(t_k-t)}{m} + \OO(n^{-1}),
\end{split} \\[0.5ex]
\begin{split}\label{rad12}
  \lambda_{jk}(1-\lambda_{jk})(1-6\lambda_{jk}+6\lambda_{jk}^2)
     &= \lambda(1-\lambda)(1-6\lambda+6\lambda^2) + \OO(n^{-1/2}).
\end{split}
\end{gather}

\nicebreak
\subsection{Estimating the factor $P(\svec,\tvec)$}\label{s:front}

Let
\[\Lambda = \prod_{j=1}^m \prod_{k=1}^n
          \lambda_{jk}^{\lambda_{jk}} (1-\lambda_{jk})^{1-\lambda_{jk}}.\]
Then
\begin{align*}
\Lambda^{-1} &= \prod_{j=1}^{m}\prod_{k=1}^n
  \biggl(\frac{1 + q_jr_k}{q_j r_k}\biggr)^{\!\!\lambda_{jk}}
        (1 + q_j r_k)^{1-\lambda_{jk}}\\
 &= \prod_{j=1}^m \prod_{k=1}^n \,(1 + q_j r_k) \,
  \biggl(\,\prod_{j=1}^m q_j^{\lambda_{j\b}}\,
       \prod_{k=1}^n r_k^{ \lambda_{\b k}}\biggr)^{\!\!-1}\\
 &= \prod_{j=1}^m \prod_{k=1}^n \,(1+q_j r_k)\,
  \prod_{j=1}^m q_j^{-s_j} \prod_{k=1}^n r_k^{-t_k}
\end{align*}
using \eqref{rad1}.  Therefore the factor
$P(\svec,\tvec)$ in front
of the integral in \eqref{start} is given by
\[ P(\svec,\tvec) = (2\pi)^{-(m+n)} \, \Lambda^{-1}.\]
We proceed to estimate $\Lambda$.
Writing $\lambda_{jk}=\lambda(1+x_{jk})$, we have
\begin{equation}\label{rad13}
\begin{split}
  &\kern-2mm\log\biggl( \frac{\lambda_{jk}^{\lambda_{jk}}
     (1-\lambda_{jk})^{1-\lambda_{jk}}}
                 {\lambda^\lambda(1-\lambda)^{1-\lambda}}\biggr)
     = \lambda x_{jk}\log\biggl(\frac{\lambda}{1-\lambda}\biggr)\\
      &\quad{} +
     \frac{\lambda}{2(1-\lambda)}x_{jk}^2 - \frac{\lambda(1-2\lambda)}
          {6(1-\lambda)^2} x_{jk}^3
       + \frac{\lambda(1-3\lambda+3\lambda^2)}{12(1-\lambda)^3} x_{jk}^4
       + O\biggl(\frac{x_{jk}^5}{(1-\lambda)^4}\biggr).
\end{split}
\end{equation}
We know from \eqref{rad1} that $\lambda_{\b\b} = mn\lambda$,
which implies that
$x_{\b\b}=0$, hence the first term on the right side of~\eqref{rad13}
does not contribute to $\Lambda$.
Now using \eqref{rad2}
we can write $x_{jk} = a_j + b_k + Z_{jk}$ and apply the
estimates in~\eqref{rad8} and~\eqref{rad9} to obtain
\begin{equation}\label{rad14}
\begin{split}
  \Lambda &= \( \lambda^\lambda(1-\lambda)^{1-\lambda} \)^{mn}
   \exp\biggr( \frac{R}{4An} + \frac{C}{4Am} + \frac{RC}{8A^2m^2n^2}
        - \frac{(1-2\lambda)R_3}{24A^2n^2} - \frac{(1-2\lambda)C_3}{24A^2m^2}\\
  &\kern 4.6cm{} + \frac{(1-3\lambda+3\lambda^2)R_4}{96A^3n^3}
               + \frac{(1-3\lambda+3\lambda^2)C_4}{96A^3m^3} + \OO(n^{-1/2})
       \biggl),
\end{split}
\end{equation}
where $R_\ell=\sum_{j=1}^m (s_j-s)^\ell$
and $C_\ell=\sum_{k=1}^n (t-t_k)^\ell$ for any~$\ell$.  Note that
$R_2=R$ and $C_2=C$.

To match the formula from the sparse case solved in~\cite{GMW}, we will
write \eqref{rad14} in terms of binomial coefficients.  First, by
Stirling's expansion of the logarithm of the gamma function,
we have that
\begin{equation}\label{rad15}
\begin{split}
   \binom{N}{(x+d)N}
  &= \frac{(x^{x+d}(1-x)^{1-x-d})^{-N}}
          {2\sqrt{\pi XN}}\\[0.5ex]
  &\quad{}\times  \exp\biggl( -\frac{1-2X}{24 XN} - \frac{d^2N}{4X}
       - \frac{(1-2x)d}{4X} + \frac{(1-4X)d^2}{16X^2} \\
      &\kern20mm {}+ \frac{(1-2x)d^3N}{24X^2}
        - \frac{(1-6X)d^4N}{96X^3}
       + O\Bigl(\frac{d^5N}{X^4} + \frac{d}{X^2N} + \frac{1}{X^3N^3}\Bigr)
  \biggr)\kern-10mm
\end{split}
\end{equation}
as $N\to\infty$, provided $x=x(N)$, $X=X(N)=\tfrac12x(1-x)$
and $d=d(N)$ are such that $0<x<1$,  $0<x+d<1$ and
provided that the error term in the above is~$o(1)$. 
{}From this we infer that
\begin{equation}\label{rad16}
\begin{split}
\binom{mn}{\lambda mn}^{\!\!-1}
  &\, \prod_{j=1}^m\binom{n}{s_j}\prod_{k=1}^n\binom{m}{t_k}
  = \frac{(\lambda^\lambda(1-\lambda)^{1-\lambda})^{-mn}}
          {(4\pi A)^{(m+n-1)/2} m^{(n-1)/2} n^{(m-1)/2}}\\
  &\quad{}\times \exp\biggl( - \frac{R}{4An} - \frac{C}{4Am}
    - \frac{1-2A}{24A}\Bigl(\frac{m}{n}+\frac{n}{m}\Bigr)
   + \frac{1-4A}{16A^2}\Bigl(\frac{R}{n^2}+\frac{C}{m^2}\Bigr)\\
  &\kern20mm{}
  + \frac{1-2\lambda}{24A^2}\Bigl(\frac{R_3}{n^2}+\frac{C_3}{m^2}\Bigr)
   - \frac{1-6A}{96A^3}\Bigl(\frac{R_4}{n^3}+\frac{C_4}{m^3}\Bigr)
   + \OO(n^{-1/2})\biggl).
\end{split}
\end{equation}
Putting \eqref{rad14} and \eqref{rad16} together, we find that
\begin{equation}\label{rad17}
\begin{split}
P(\svec,\tvec)
   &= \Lambda^{-1}(2\pi)^{-(m+n)}\\
   &= \frac{A^{(m+n-1)/2} m^{(n-1)/2} n^{(m-1)/2}}
                       {2\pi^{(m+n+1)/2}}
      \;\binom{mn}{\lambda mn}^{\!\!-1}
      \prod_{j=1}^{m}\binom{n}{s_j}\prod_{k=1}^{n}\binom{m}{t_k}  \\
   & \hspace*{5mm}{} \times
     \exp\biggl( \frac{1-2A}{24A}\Bigl(\frac{m}{n}+\frac{n}{m}\Bigr)
      - \frac{RC}{8A^2m^2n^2}
      - \frac{1-4A}{16A^2}\Bigl(\frac{R}{n^2}+\frac{C}{m^2}\Bigr)
      + \OO(n^{-1/2})\biggr).\kern-3mm
\end{split}
\end{equation}

\nicebreak
\section{Evaluating the integral}\label{s:integral}

Our next task is to evaluate the integral $I(\svec,\tvec)$
given by
\begin{equation}
I(\svec,\tvec) = \int_{-\pi}^\pi \!\!\cdots \int_{-\pi}^\pi
\frac{ \prod_{j=1}^m\prod_{k=1}^n 
\((1 + \lambda_{j,k}(e^{i(\t_j+\p_k)} - 1) \)}
{\exp(i\sum_{j=1}^m s_j\theta_j + i\sum_{k=1}^n t_k\phi_k)}\,
    d \thetavec  d\phivec.
\label{FDef}
\end{equation}

It is convenient to think of $\theta_j$, $\phi_k$ as points on
the unit circle.
We wish to define ``averages'' of the angles $\theta_j$, $\phi_k$.
To do this cleanly we make the following definitions, as in~\cite{CM}.
Let $C$ be the ring of real numbers modulo $2\pi$, which we can interpret
as points on a circle in the usual way.  Let $z$ be the canonical
mapping from $C$ to the real interval $(-\pi, \pi]$.  An \emph{open
half-circle} is $C_t = (t-\pi/2, t+\pi/2) \subseteq C$ for some $t$.
Now define
\[ \Chat^N = \bigcup_t \Chat_t^N = 
  \{\, \xvec = (x_1,\ldots, x_N) \in C^N \mid
  x_1, \ldots, x_N \in C_t \mbox{ for some } t\in\mathbb{R}\,\}.\]
If $\xvec = (x_1,\ldots, x_N)\in C_0^N$ then define
\[ \avx = z^{-1}\biggl( \frac{1}{N}
  \sum_{j=1}^N z(x_j)\biggr).\]
More generally, if $\xvec\in C_t^N$ then define
$\avx = t + \overline{(x_1-t,\ldots, x_N-t)}$.
The function $\xvec\to \avx$ is
well-defined and continuous for $\xvec\in\Chat^N$.

Let $\R$
denote the set of vector pairs $(\thetavec ,\phivec)
\in \Chat^m\times\Chat^n$ such that
\begin{align}\label{Rdef}
\begin{split}
\abs{\avtheta + \avphi}
  &\leq  (mn)^{-1/2 + 2\eps},\\
\abs{\that_j} &\leq n^{-1/2 + \eps} \qquad (1\leq j\leq m),\\
\abs{\phat_k} &\leq m^{-1/2 + \eps} \qquad (1\leq k\leq n),
\end{split}
\end{align}
where $\that_j = \theta_j - \avtheta$
and $\phat_k = \phi_k - \avphi$.  In
this definition, values are considered in $C$.  The constant
$\eps$ is the sufficiently-small value required by
Theorem~\ref{bigtheorem}.

\smallskip

Let $I_{\R''}(\svec,\tvec)$ denote the integral
$I(\svec,\tvec)$ restricted to any region $\R''$.
In this section, we estimate $I_{\R'}(\svec,\tvec)$
in a certain region $\R'\supseteq\R$.
In Section~\ref{s:boxing} we will show that
the remaining parts of $I(\svec,\tvec)$ are negligible.
We begin by analysing
the integrand in $\R$, but for future use when we
expand the region to~$\R'$ (to be defined in~\eqref{RSprime}),
note that all the approximations
we establish for the integrand in $\R$ also hold in the
superset of $\R'$ defined~by
\begin{align}\begin{split}\label{bigR}
\abs{\avtheta + \avphi}
  &\leq  3(mn)^{-1/2 + 2\eps},\\
\abs{\that_j} &\leq 3n^{-1/2 + \eps} \quad (1\leq j\leq m-1),\\
\abs{\that_m} &\leq 2n^{-1/2 + 3\eps},\\
\abs{\phat_k} &\leq 3m^{-1/2 + \eps} \quad (1\leq k\leq n-1),\\
\abs{\phat_n} &\leq 2m^{-1/2 + 3\eps}.
\end{split}\end{align}

Define $\thetahatvec=(\that_1,\ldots,\that_{m-1})$ and
$\phihatvec=(\phat_1,\ldots,\phat_{n-1})$.
Let $T_1$ be the transformation
$T_1(\thetahatvec,\phihatvec,\nu,\delta)=(\thetavec,\phivec)$
defined by
\[ \nu = \avtheta + \avphi,
  \quad  \delta = \avtheta -
\avphi,\]
together with $\that_j=\t_j-\avtheta$
($1\leq j\leq m-1$) and $\phat_k=\p_k-\avphi$ ($1\leq k\leq n-1$).
We also define the 1-many transformation $T_1^*$ by
\[
    T_1^*(\thetahatvec,\phihatvec,\nu) =
    \bigcup_\delta \, T_1(\thetahatvec,\phihatvec,\nu,\delta).
\]

After applying the transformation $T_1$ to $I_\R(\svec ,\tvec)$,
the new integrand is easily seen to be independent of
$\delta$, so we can multiply by the range of $\delta$ and remove
it as an independent variable.  Therefore, we can continue with
an $(m+n-1)$-dimensional integral over $\S$ such that
$\R=T_1^*(\S)$.   More generally, if
$\S''\subseteq (-\tfrac12\pi,\tfrac12\pi)^{m+n-2}\times(-2\pi,2\pi]$
and $\R''=T_1^*(\S'')$, we have
\begin{equation}\label{IvsJ}
I_{\R''}(\svec,\tvec) =  2\pi m n \int_{\S''} G(\thetahatvec,\phihatvec,\nu)
     \, d\thetahatvec  d\phihatvec d\nu,
\end{equation}
where $G(\thetahatvec,\phihatvec,\nu)
=F\(T_1(\thetahatvec,\phihatvec,\nu,0)\)$
with $F(\thetavec,\phivec)$ being the integrand of~\eqref{FDef}.
The factor $2\pi mn$ combines the range of $\delta$,
which is $4\pi$, and the Jacobian of $T_1$, which is~$mn/2$.

Note that  $\S$ is defined by
the same inequalities \eqref{Rdef} as define~$\R$.  The first inequality
is now $\abs{\nu}\leq (mn)^{-1/2 + 2\eps}$ and the bounds on
\[ \that_m = - \jsum \that_j~\text{ and }~
    \phat_m = - \ksum \phat_k
\]
still apply even though these are no
longer variables of integration.

\medskip

Our main result in this section is the following.
\begin{thm}\label{Jintegral}
Under the conditions of Theorem~\ref{bigtheorem},  there is a
region $\S'\supseteq\S$ such that
\begin{align*}
 \int_{\S'} G(\thetahatvec,\phihatvec,\nu)
     \, d\thetahatvec  d\phihatvec  d\nu
     &= (mn)^{-1/2} \Bigl(\frac{\pi}{Amn }\Bigr)^{\!1/2}
                    \Bigl(\frac{\pi}{An  }\Bigr)^{(m-1)/2}
                    \Bigl(\frac{\pi}{Am  }\Bigr)^{(n-1)/2} \\
     &\kern-18mm {}\times \exp\biggl(
         - \frac12
         - \frac{1-2A}{24A}\Bigl(\frac mn+\frac nm\Bigr)
   + \frac{1}{4A}\Bigl(\frac{1}{m}+\frac{1}{n}\Bigr)
               \Bigl(\frac{R}{n}+\frac{C}{m}\Bigr) \\
      &\kern32mm{}
      + \frac{1-8A}{16A^2}\Bigl(\frac{R}{n^2}+\frac{C}{m^2}\Bigr)
       + O(n^{-b})  \biggr).
\end{align*}
\end{thm}

\medskip

In the region $\S$, the integrand of \eqref{IvsJ} can be expanded as
\begin{align*}
G(\thetahatvec,\phihatvec,\nu) &= \exp\biggl( -\sumjkmn
       (A+\alpha_{jk})(\nu+\that_j+\phat_k)^2
                 - i\sumjkmn(A_3+\beta_{jk})(\nu+\that_j+\phat_k)^3
\\
    & \qquad\qquad{}
                 + \sumjkmn(A_4+\gamma_{jk})(\nu+\that_j+\phat_k)^4
                 + O\Bigl(A\,\sumjkmn\,\abs{\nu+\that_j+\phat_k}^5\,\Bigr)
             \biggr).
\end{align*}
Here $\alpha_{jk}$, $\beta_{jk}$, and $\gamma_{jk}$ are defined by
\begin{align}\label{AlBetGamDef}
\begin{split}
\dfrac{1}{2}\lambda_{jk}(1-\lambda_{jk}) &=  A + \alpha_{jk}, \\
\dfrac{1}{6}\lambda_{jk}(1-\lambda_{jk})(1-2\lambda_{jk})
   &= A_3 + \beta_{jk}, \\
\dfrac{1}{24}\lambda_{jk}(1-\lambda_{jk})(1-6\lambda_{jk}+6\lambda_{jk}^2)
   &= A_4 + \gamma_{jk},
\end{split}
\end{align}
where
\[ A = \dfrac{1}{2}\lambda(1-\lambda),~
 A_3 = \dfrac{1}{6}\lambda(1-\lambda)(1-2\lambda), \text{~and~}
  A_4 = \dfrac{1}{24}\lambda(1-\lambda)(1-6\lambda + 6\lambda^2).
\]
Approximations for $\alpha_{jk}$, $\beta_{jk}$, $\gamma_{jk}$ were
given in \eqref{rad10}--\eqref{rad12}.

\nicebreak
\subsection{Another change of variables}\label{s:change}

We now make a second change of variables
$(\thetahatvec,\phihatvec,\nu) =T_2(\zetavec,\xivec,\nu)$,
where $\zetavec=(\zeta_1,\ldots,\zeta_{m-1})$ and
$\xivec=(\xi_1,\ldots,\xi_{n-1})$,
whose purpose is to almost diagonalize the quadratic part
of~$G$.  The diagonalization will be completed in the
next subsection.
The transformation $T_2$ is defined as follows.
For $1\leq j\leq m-1$ and $1\leq k\leq n-1$ let
\[ \that_j = \zeta_j + c\pi_1,\quad \phat_k = \xi_k + d\rho_1, \]
where
\[ c=-\frac{1}{m+m^{1/2}}
   \quad \mbox{and} \quad d= - \frac{1}{n+n^{1/2}}\]
and, for $1\leq h\leq 4$,
\[ \pi_h = \jsum \zeta_j^h,\quad \rho_h = \ksum
           \xi_k^h.\]
The Jacobian of the transformation is $(mn)^{-1/2}$.
In \cite{CM}, this transformation was seen to exactly diagonalize
the quadratic part of the integrand in the semiregular case.
In the present irregular case, the diagonalization is no longer
exact but still provides useful progress.

By summing the equations $\that_j = \zeta_j + c\pi_1$ and
$\phat_k = \xi_k + d\rho_1$, we find that
\begin{align}\label{pi1bound}
\begin{split}
   \pi_1 = m^{1/2}\jsum\that_j,\quad \abs{\pi_1}\le m^{1/2}n^{-1/2+\eps}, \\
   \rho_1 = n^{1/2}\ksum\phat_k,\quad \abs{\rho_1}\le n^{1/2}m^{-1/2+\eps},
\end{split}
\end{align}
where the right sides come from the bounds on $\that_m$ and $\phat_n$.
This implies that
\begin{align}\label{zetaxi}
\begin{split}
   \zeta_j &= \that_j+\OO(n^{-1})\quad(1\le j\le m-1), \\
   \xi_k &= \phat_k+\OO(n^{-1})\quad(1\le k\le n-1).
\end{split}
\end{align}
The transformed region of integration is $T_2^{-1}(\S)$, but
for convenience we will expand it a little to be the region
defined by the inequalities
\begin{align}\label{zetaxibox}
\begin{split}
   \abs{\zeta_j} &\le \tfrac32 n^{-1/2+\eps}\quad(1\le j\le m-1), \\
   \abs{\xi_k} &\le \tfrac32 m^{-1/2+\eps}\quad(1\le k\le n-1), \\
   \abs{\pi_1} &\le m^{1/2}n^{-1/2+\eps}, \\
   \abs{\rho_1} &\le n^{1/2}m^{-1/2+\eps}, \\
   \abs{\nu} &\le (mn)^{-1/2+2\eps}\,.
\end{split}
\end{align}

We now consider the new integrand $E_1=\exp(L_1) = G\circ T_2$.
As in~\cite{CM}, the semiregular parts of the integrand (those not
involving $\alpha_{jk}$, $\beta_{jk}$ or~$\gamma_{jk}$) transform to 
\begin{align}\label{regbits}
\begin{split}
&-Amn\nu^2 - A n \pi_2 - A m \rho_2
   -3i A_3 n\nu \pi_2 - 3iA_3 m\nu \rho_2
  + 6 A_4\pi_2\rho_2 \\
 &\quad  {} - i A_3 n\pi_3
        - i A_3 n\rho_3
- 3i A_3 cn \pi_1\pi_2
- 3iA_3 dm \rho_1\rho_2
 + A_4 n\pi_4 + A_4 m\rho_4 + \OO(n^{-1/2}).
\end{split}
\end{align}
To see the effect of the transformation on the irregular parts of the
integrand, write $\zeta_m=\that_m-c\pi_1$ and
$\xi_n=\that_n-d\rho_1$.  {}From~\eqref{pi1bound} we can see
that $\zeta_m=\OO(n^{-1/2})$ and $\xi_n=\OO(n^{-1/2})$.
Thus we have, for all $1\le j\le m$ and $1\le k\le n$,
$\zeta_j+\xi_k = \OO(n^{-1/2})$ and $c\pi_1+d\rho_1+\nu=\OO(n^{-1})$.
Recalling also that $\alpha_{jk}, \beta_{jk}, \gamma_{jk} = \OO(n^{-1/2})$,
we have
\begin{align*}
  \sum_{j=1}^m\sum_{k=1}^n\alpha_{jk} (\nu+\that_j+\phat_k)^2\\[-2ex]
  &\kern-2cm{}= \sum_{j=1}^m\sum_{k=1}^n\alpha_{jk} 
         \((\zeta_j+\xi_k)^2+2(\zeta_j+\xi_k) (\nu+c\pi_1+d\rho_1)\) + \OO(n^{-1/2}), \\
  \sum_{j=1}^m\sum_{k=1}^n\beta_{jk} (\nu+\that_j+\phat_k)^3
  &= \sum_{j=1}^m\sum_{k=1}^n\beta_{jk}(\zeta_j+\xi_k)^3 + \OO(n^{-1/2}), \\
  \sum_{j=1}^m\sum_{k=1}^n\gamma_{jk} (\nu+\that_j+\phat_k)^4
  &= \OO(n^{-1/2}).
\end{align*}
Moreover, the terms on the right sides of the above that involve
$\zeta_m$ or $\xi_n$ contribute only $\OO(n^{-1/2})$ in total, so we can drop them.
Combining this with~\eqref{regbits}, we have
\begin{equation}\label{Helper}
\begin{split}
L_1 &=  -Amn\nu^2 - A n \pi_2 - A m \rho_2
   -3i A_3 n\nu \pi_2 - 3iA_3 m\nu \rho_2
  + 6 A_4\pi_2\rho_2 \\[0.6ex]
 & \hspace*{5mm} {} - i A_3 n\pi_3
        - i A_3 n\rho_3
- 3i A_3 cn \pi_1\pi_2
- 3iA_3 dm \rho_1\rho_2
 + A_4 n\pi_4 + A_4 m\rho_4 \\
 & \hspace*{5mm} {}
   - \jsum\ksum
   \alpha_{jk}\((\zeta_j+\xi_k)^2+2(\zeta_j+\xi_k) (\nu+c\pi_1+d\rho_1)\)
          \\
 & \hspace*{5mm} {}
   - i \jsum\ksum \beta_{jk}(\zeta_j+\xi_k)^3
   + \OO(n^{-1/2}).
\end{split}
\end{equation}

\nicebreak
\subsection{Completing the diagonalization}\label{s:diagonalize}

The quadratic form in $E_1$ is the following function of the $m+n-1$
variables $\zetavec,\xivec,\nu$:
\begin{align}\begin{split}\label{Qvalue}
Q &= -Amn\nu^2 - An \pi_2 - Am \rho_2 \\
    &\quad {} - \jsum\ksum
\alpha_{jk}\((\zeta_j+\xi_k)^2 + 2(\zeta_j+\xi_k)(\nu+c\pi_1+d\rho_1)\).
\end{split}\end{align}
We will make a third change of variables,
$(\zetavec,\xivec,\nu)=T_3 (\sigmavec,\tauvec,\mu)$, that
diagonalizes this quadratic form,
where $\sigmavec=(\sigma_1,\ldots,\sigma_{m-1})$
and $\tauvec=(\tau_1,\ldots,\tau_{n-1})$.
This is achieved using a slight extension
of \cite[Lemma 3.2]{MWtournament}.

\begin{lemma}
\label{diagonalize}
Let $X$ and $Y$ be square matrices of the same order, such that
$X^{-1}$ exists and all the eigenvalues of
$X^{-1} Y$ are less than 1 in absolute value.
Then
\[ (I + YX^{-1})^{-1/2} \,(X+Y)\, (I+X^{-1}Y)^{-1/2} = X, \]
where the fractional powers are defined by the binomial expansion.
\quad\qedsymbol
\end{lemma}

Note that $X^{-1} Y$ and $YX^{-1}$ have the same eigenvalues, so
the eigenvalue condition on $X^{-1} Y$ applies equally to $YX^{-1}$.
If we also have that both $X$ and $Y$ are symmetric, then
\[ \sum_{r\geq 0}\binom{-\dfrac12}{\,r} (YX^{-1})^T
      = \sum_{r\geq 0} \binom{-\dfrac12}{\,r} (X^{-1})^T Y^T
         = \sum_{r\geq 0} \binom{-\dfrac12}{\,r} X^{-1} Y \]
so $(I+YX^{-1})^{-1/2}$ is the transpose of $(I + X^{-1}Y)^{-1/2}$.
Let $V$ be the symmetric matrix associated with the quadratic form
$Q$.   Write $V=\Vd + \Vnd$ where $\Vd$ has all
off-diagonal
entries equal to zero and matches $V$ on the diagonal entries,
and $\Vnd$ has all diagonal entries zero and matches $V$ on the
off-diagonal entries.
We will apply Lemma~\ref{diagonalize}
with $X = \Vd$ and $Y=\Vnd$.
Note that $\Vd$ is invertible and that both $\Vd$
and $\Vnd$ are symmetric.
Let $T_3$ be the transformation given by
$T_3(\sigmavec,\tauvec,\mu)^T = (\zetavec,\xivec,\nu)^T=
(I+\Vdinv \Vnd)^{-1/2}(\sigmavec,\tauvec,\mu)^T$.
If the eigenvalue condition of Lemma~\ref{diagonalize} is satisfied then
this transformation diagonalizes the quadratic form $Q$, keeping
the diagonal entries unchanged.

{}From the formula for $Q$ we extract the following coefficients,
which tell us the diagonal and off-diagonal entries of $V$:
\begin{align*}
[\zeta_j^2]\,Q &= -An-(1+2c)\alpha_{j\c},\\
[\xi_k^2]\,Q &= -Am-(1+2d)\alpha_{\c k},\\
[\nu^2]\,Q &= -Amn, 
\displaybreak[0]\\
[\zeta_{j_1}\zeta_{j_2}]\,Q &= -2c(\alpha_{j_1\c}+\alpha_{j_2\c})
 \qquad(j_1\ne j_2),\\
[\zeta_j \xi_k]\,Q &= -2\alpha_{jk} - 2d\alpha_{j\c} - 2c\alpha_{\c k},\\
[\xi_{k_1}\xi_{k_2}]\,Q &= -2d(\alpha_{\c k_1}+\alpha_{\c k_2})
 \qquad(k_1\ne k_2),\\
[\zeta_j \nu]\,Q &= -2\alpha_{j\c}, \\
[\xi_k \nu]\,Q &= -2\alpha_{\c k}.
\end{align*}
Using these equations we find that all off-diagonal entries of
$\Vdinv\Vnd$ are $\OO(n^{-3/2})$, except for the column
corresponding to $\nu$ which has off-diagonal entries of size $\OO(n^{-1/2})$.
Similarly, the off-diagonal entries of $\Vnd\Vdinv$
are all $\OO(n^{-3/2})$, except for the row corresponding to $\nu$,
which has off-diagonal entries of size $\OO(n^{-1/2})$.  To see that
these conditions imply that the eigenvalues of
$\Vdinv\Vnd$
are less than one, recall that the value of any matrix norm is
greater than or equal to the greatest absolute value of an
eigenvalue.  The $\infty$-norm (maximum row sum of
absolute values) of $\Vdinv\Vnd$ is $\OO(n^{-1/2})$, so the
eigenvalues are all~$\OO(n^{-1/2})$.

We also need to know the Jacobian of the transformation $T_3$.

\begin{lemma}\label{l:detexp}
  Let $M$ be a matrix of order $O(m+n)$ with all eigenvalues
  uniformly $\OO(n^{-1/2})$.  Then
  \[
    \det(I+M) = \exp\(\tr M - \dfrac12\tr M^2 + \OO(n^{-1/2})\).
  \]
\end{lemma}
\proof
The eigenvalue condition ensures that the Taylor series for
$\log(I+M)$ converges and that
\[
 \det(I+M) = \exp\(\tr\log(I+M)\).
\]
Expanding the logarithm and noting that
$\abs{\tr M^r}=\OO(n^{-(r-2)/2})$ for $r\ge 3$ gives the result.
\endproof

Let $M=\Vdinv\Vnd$.
As noted before, the eigenvalues of $M$ are all $\OO(n^{-1/2})$ so
Lemma~\ref{l:detexp} applies.
Noting that
$\tr(M)=0$ and calculating that $\tr(M^2) = \OO(n^{-1})$, we
conclude that the Jacobian of $T_3$ is
\[ \det\((I+M)^{-1/2}\) = \(\det(I+M)\)^{-1/2} = 1+ \OO(n^{-1/2}).  \]

\medskip

To derive $T_3$ explicitly, we can expand
$(I+\Vdinv\Vnd)^{-1/2}$ while noting that
$\alpha_{j\c} = O(n^{1/2+\eps})$ for all~$j$,
$\alpha_{\c k} = O(m^{1/2 + \eps})$ for all~$k$,
$\alpha_{\c\c} = O(mn^{2\eps} + nm^{2\eps})$,
$R\le mn^{1+2\eps}$ and $C\le nm^{1+2\eps}$.

\nicebreak
This gives
\begin{align*}
\sigma_j &= \zeta_j + \sum_{j'=1}^{m-1} \Bigl(\frac{c(\alpha_{j\c} +
        \alpha_{j'\c})}{2An} + \OO(n^{-2})\Bigr) \zeta_{j'} +
        \ksum \Bigl(\frac{\alpha_{jk} + d\alpha_{j\c} +
          c\alpha_{\c k}}{2An} + \OO(n^{-2})\Bigr) \xi_k \\
      & \hspace*{2cm} {} +
         \Bigl(\frac{\alpha_{j\c}}{2An} + \OO(n^{-1})\Bigr)\nu + \OO(n^{-2}),
       \displaybreak[0]\\
\tau_k &= \xi_k + \jsum \Bigl(\frac{\alpha_{jk} +
      d\alpha_{j\c} + c\alpha_{\c k}}{2Am} + \OO(n^{-2})\Bigr)\zeta_j
    + \sum_{k'=1}^{n-1} \Bigl(\frac{d(\alpha_{\c k} + \alpha_{\c k'})}
  {2Am} + \OO(n^{-2})\Bigr)\xi_{k'} \\
      & \hspace*{2cm} {} + \Bigl(\frac{\alpha_{\c k}}{2Am} + \OO(n^{-1})
         \Bigr)\nu + \OO(n^{-2}),
         \displaybreak[0]\\
\mu &= \nu + \jsum
   \Bigl(\frac{\alpha_{j\c}}{2Amn} + \OO(n^{-2})\Bigr)\zeta_j
     + \ksum \Bigl( \frac{\alpha_{\c k}}{2Amn} +\OO(n^{-2})\Bigr)\xi_k
        + \OO(n^{-1})\nu,
\end{align*}
for $1\leq j\leq m-1$, $1\leq k\leq n-1$.

\medskip

The transformation $T_3^{-1}$ perturbs the region of integration in
an irregular  fashion that we must bound.  {}From the explicit form of
$T_3$ above, we have
\begin{align*}
\sigma_j &= \zeta_j + \sum_{j' = 1}^{m-1} \OO(n^{-3/2}) \zeta_{j'} +
                          \ksum \OO(n^{-3/2}) \xi_k + \OO(n^{-1/2}) \nu+ \OO(n^{-2})
         = \zeta_j + \OO(n^{-1}),\\
\tau_k &= \xi_k + \sum_{j = 1}^{m-1} \OO(n^{-3/2}) \zeta_{j} +
                          \sum_{k'=1}^{n-1} \OO(n^{-3/2}) \xi_{k'} 
                          + \OO(n^{-1/2}) \nu+ \OO(n^{-2})
       = \xi_k + \OO(n^{-1})
\end{align*}
for $1\leq j\leq m-1$, $1\leq k\leq n-1$,
so $\sigmavec,\tauvec$ are only slightly different from $\zetavec,\xivec$.

For $\mu$ versus $\nu$ we have
\begin{align*}
  \mu &= \nu  + O(n^{-1+2\eps}/A) + O(m^{-1+2\eps}/A) \\
          &= \nu + o\( (mn)^{-1/2+2\eps}\),
\end{align*}
where the second step requires our assumptions
$m=o(A^2n^{1+\eps})$ and $n= o(A^2 m^{1+\eps})$.
This shows that the bound $\abs{\nu}\le  (mn)^{-1/2+2\eps}$
is adequately covered by
$\abs\mu \le 2(mn)^{-1/2 + 2\eps}$.

For $1\leq h\leq 4$, define
\[ \mu_h = \jsum {\sigma_j}^h,\quad \nu_h = \ksum
            {\tau_k}^h.\]
{}From \eqref{zetaxibox}, we see that
$\abs{\pi_1} \le m^{1/2} n^{-1/2+ \eps}$ and
$\abs{\rho_1} \le m^{-1/2 + \eps}n^{1/2}$ are the remaining
constraints that define the region of integration.  We next apply
these constraints to bound $\mu_1$ and~$\nu_1$.
{}From the explicit form of $T_3$, we have
\begin{align}
 \mu_1 &= \pi_1
  + \jsum \sum_{j'=1}^{m-1}
  \Bigl(\frac{c(\alpha_{j\c}+
         \alpha_{j'\c})}{2An}+ \OO(n^{-2})\Bigr)\zeta_{j'}\notag 
    \\
  &  \kern10mm {}
   + \jsum\ksum
  \Bigl(\frac{\alpha_{jk} + d\alpha_{j\c} + c\alpha_{\c k}}
                      {2An} + \OO(n^{-2})\Bigr) \xi_k
  + \jsum \Bigl(\frac{\alpha_{j\c}}{2An}
  + \OO(n^{-1})\Bigr)\nu + \OO(n^{-1}) \notag \displaybreak[0]\\
 &= \pi_1 + \frac{c\alpha_{\c\c}}{2An} m^{1/2} n^{-1/2+\eps}
 + \frac{d\alpha_{\c\c}}{2An} m^{-1/2+\eps}n^{1/2}
      + \frac{\alpha_{\c\c}}{2An} \nu \notag \\
& \kern10mm  + \(1+c(m-1)\)\ksum\frac{\alpha_{\c k}}{2An}\xi_k
 + \frac{c(m-1)}{2An}\sum_{j'=1}^{m-1} \alpha_{j'\c}\zeta_{j'}
  + \OO(n^{-1/2}) \notag \displaybreak[0]\\
 &= \pi_1
  + \frac{c(m-1)}{2An} \sum_{j'=1}^{m-1} \alpha_{j'\c}\zeta_{j'}
  + \OO(n^{-1/2}) \label{mu1pi1}\\
  &= \pi_1 +O(A^{-1}mn^{-1+2\eps})\notag  \\[0.5ex]
  &= \pi_1 + o(m^{1/2}n^{-1/2+5\eps/2}).\notag
\end{align}
To derive the above we have used $1+c(m-1)=m^{1/2}$ and the bounds we
have established on the various variables.  For the last step, we
need the assumption $m=o(A^2n^{1+\eps})$, which implies that
$A^{-1}mn^{-1+2\eps}=o(m^{1/2}n^{-1/2+5\eps/2})$.

Since our region of integration has
$\abs{\pi_1}\le m^{1/2}n^{-1/2+\eps}$, we see that this implies
the bound $\abs{\mu_1}\le m^{1/2}n^{-1/2+3\eps}$.
By a parallel argument, we have
\[\nu_1 = \rho_1 + o(m^{-1/2+5\eps/2}n^{1/2}),\]
which implies
$\abs{\nu_1}\le n^{1/2}m^{-1/2+3\eps}$.
Putting together all the bounds we have derived, we see that
\[T_3^{-1}(T_2^{-1}(\S)) \subseteq \Q\cap\M,\]
where
\begin{align*}
  \Q &= \{\, \abs{\sigma_j}\le 2n^{-1/2 + \eps}, j=1,\ldots, m-1\,\}
   \cap \{\, \abs{\tau_k}\le 2m^{-1/2 + \eps}, k=1,\ldots, n-1\,\}\\
  & \hspace*{1cm} \cap \{ \abs\mu \le 2(mn)^{-1/2+2\eps}\,\},\\
  \M &= \{\, \abs{\mu_1}\le m^{1/2} n^{-1/2+3\eps}\,\} \cap
        \{\, \abs{\nu_1}\le n^{1/2} m^{-1/2+3\eps}\}.
\end{align*}

Now define
\begin{align}\label{RSprime}
\begin{split}
 \S' &=T_2(T_3(\Q\cap\M)), \\
 \R' &= T_1^*(\S').
\end{split}
\end{align}
 We have proved that
$\S'\supseteq\S$, so it is valid to take~$\S'$
to be the region required by Theorem~\ref{Jintegral}.
Also notice that $\R'$ is contained in the region
defined by the inequalities~\eqref{bigR}.
As we forecast at that time, our estimates of the integrand
have been valid inside this expanded region.  It remains
to apply the transformation $T_3^{-1}$ to the
integrand~\eqref{Helper} so that we have it in terms of
$(\sigmavec,\tauvec,\mu)$.  The explicit form of $T_3^{-1}$
is similar to the explicit form for $T_3$, namely:
\begin{align*}
\zeta_j &= \sigma_j - \sum_{j'=1}^{m-1} \Bigl(\frac{c(\alpha_{j\c} +
        \alpha_{j'\c})}{2An} + \OO(n^{-2})\Bigr) \sigma_{j'} -
        \ksum \Bigl(\frac{\alpha_{jk} + d\alpha_{j\c} +
          c\alpha_{\c k}}{2An} + \OO(n^{-2})\Bigr) \tau_k \\
      & \hspace*{2cm} {} -
         \Bigl(\frac{\alpha_{j\c}}{2An} + \OO(n^{-1})\Bigr)\mu +\OO(n^{-2}),
       \displaybreak[0]\\[0.5ex]
\xi_k &= \tau_k - \jsum \Bigl(\frac{\alpha_{jk} -
      d\alpha_{j\c} + c\alpha_{\c k}}{2Am} + \OO(n^{-2})\Bigr)\sigma_j
    - \sum_{k'=1}^{n-1} \Bigl(\frac{d(\alpha_{\c k} + \alpha_{\c k'})}
  {2Am} + \OO(n^{-2})\Bigr)\tau_{k'} \\
      & \hspace*{2cm} {} - \Bigl(\frac{\alpha_{\c k}}{2Am} + \OO(n^{-1})
         \Bigr)\mu +\OO(n^{-2}),\\[0.5ex]
\nu &= \mu - \jsum
   \Bigl(\frac{\alpha_{j\c}}{2Amn} + \OO(n^{-2})\Bigr)\sigma_j
     - \ksum \Bigl( \frac{\alpha_{\c k}}{2Amn} +\OO(n^{-2})\Bigr)\tau_k
        + \OO(n^{-1})\mu,
\end{align*}
for $1\leq j\leq m-1$, $1\leq k\leq n-1$.
In addition to the relationships between the old and new variables
that we proved before, we can note that $\pi_2=\mu_2+\OO(n^{-1/2})$,
$\rho_2=\nu_2+\OO(n^{-1/2})$, $\pi_3=\mu_3+\OO(n^{-1})$,
$\rho_3=\nu_3+\OO(n^{-1})$,  $\pi_4=\mu_4+\OO(n^{-3/2})$,
and $\rho_4=\nu_4+\OO(n^{-3/2})$.

\medskip
The quadratic part of $L_1$, which we called $Q$ in~\eqref{Qvalue},
loses its off-diagonal parts according to our design of $T_3$.
Thus, what remains is
\begin{align*}
 -Amn\mu^2 &- \jsum \(An+(1+2c)\alpha_{j\c}\)\sigma_j^2
   - \ksum \(Am+(1+2d)\alpha_{\c k}\)\tau_k^2\\
     &= -Amn\mu^2 - An\mu_2 - Am\nu_2
         - \jsum \alpha_{j\c} \sigma_j^2
         - \ksum \alpha_{\c k} \tau_k^2 + \OO(n^{-1/2}).
\end{align*}

Next consider the cubic terms of $L_1$.  These are
\begin{align*}
   &{}-3i A_3 n\nu \pi_2 - 3iA_3 m\nu \rho_2
   - i A_3 n\pi_3 - i A_3 n\rho_3 \\
&\qquad{}- 3i A_3 cn \pi_1\pi_2 - 3iA_3 d n \rho_1\rho_2 
- i \jsum\ksum \beta_{jk}(\zeta_j+\xi_k)^3.
\end{align*}
We calculate the following in $\Q\cap\M$:
\begin{align}
 -3i A_3 n\nu \pi_2 &= -3i A_3 n\mu\mu_2
       + \frac{3iA_3\mu_2}{2Am}\biggl(\,\jsum \alpha_{j\c}\sigma_j
                    + \ksum \alpha_{\c k}\tau_k\biggr) + \OO(n^{-1/2}), \notag\\
  -i A_3 n\pi_3 &= -i A_3 n\mu_3 + \frac{3iA_3}{2A}\biggl(\,
       \sum_{j,j'=1}^{m-1}
         c(\alpha_{j\c}+\alpha_{j'\c})\sigma_j^2\sigma_{j'}, \notag \\
    &\kern40mm{}   + \jsum\ksum
           (\alpha_{jk}+d\alpha_{j\c}+c\alpha_{\c k})
             \sigma_j^2\tau_k\biggr) + \OO(n^{-1/2}), \notag 
             \displaybreak[0]\\
  - 3i A_3 cn \pi_1\pi_2 &= - 3i A_3 cn \mu_1\mu_2
     + \frac{3iA_3c^2m\mu_2}{2A}
           \jsum \alpha_{j\c}\sigma_j + \OO(n^{-1/2}) \label{pi1pi2}, \\
  - i \jsum\ksum \beta_{jk}&(\zeta_j+\xi_k)^3 =
     - i \jsum\ksum \beta_{jk}(\sigma_j+\tau_k)^3
         + \OO(n^{-1/2}), \label{beta3sum}
\end{align}
and the remaining cubic terms are each parallel to one of those.
The proof of \eqref{pi1pi2} is similar to the proof of \eqref{mu1pi1}.

Finally we come to the quartic part of $E_1$, which is
\[6A_4\pi_2\rho_2 + A_4n\pi_4 + A_4m\rho_4
  = 6A_4\mu_2\nu_2 + A_4n\mu_4 + A_4m\nu_4+ \OO(n^{-1/2}).\]

In summary, the value of the integrand for
$(\sigmavec,\tauvec,\mu)\in\Q\cap\M$ is
$\exp\(L_2+\OO(n^{-1/2})\)$, where
\begin{align}\label{L2value}
\begin{split}
   L_2 &= -Amn\mu^2 - An\mu_2 - Am\nu_2
         - \jsum \alpha_{j\c} \sigma_j^2
         - \ksum \alpha_{\c k} \tau_k^2
         + 6A_4\mu_2\nu_2 \\
         &\quad{}+ A_4n\mu_4 + A_4m\nu_4
      - i A_3 n\mu_3 - i A_3 m\nu_3
               - 3i A_3 cn \mu_1\mu_2 - 3i A_3 dm \nu_1\nu_2\\
       &\quad{} - 3 i A_3 n \mu\mu_2 - 3 i A_3 m \mu\nu_2
                - i \jsum \beta_{j\c}\sigma_j^3
               - i \ksum \beta_{\c k}\tau_k^3 \\
       &\quad{} + i\sum_{j,j'=1}^{m-1} g_{jj'}\sigma_j\sigma_{j'}^2
             + i\sum_{k,k'=1}^{n-1} h_{kk'}\tau_k\tau_{k'}^2
             + i\jsum\ksum
                \(u_{jk}\sigma_j\tau_k^2 + v_{jk}\sigma_j^2\tau_k\),
\end{split}
\end{align}
with
\begin{align*}
    g_{jj'} &= \frac{3A_3}{2Am}
                    \((1+cm+c^2m^2)\alpha_{j\c}+cm\alpha_{j'\c}\)
                    = O(n^{-1/2+\eps}), \\[0.4ex]
    h_{kk'} &= \frac{3A_3}{2An}
                    \((1+dn+d^2n^2)\alpha_{\c k}+dn\alpha_{\c k'}\)
                    = O(m^{-1/2+\eps}), \\[0.4ex]
    u_{jk} &= \frac{3A_3}{2An}
                    \( n\alpha_{jk} + (1+dn)\alpha_{j\c} + cn\alpha_{\c k}\)
                   - 3\beta_{jk} = O(m^{-1/2+2\eps}+n^{-1/2+2\eps}), \\[0.4ex]
    v_{jk} &= \frac{3A_3}{2Am}
                    \( m\alpha_{jk} + (1+cm)\alpha_{\c k} + dm\alpha_{j\c}\)
                           - 3\beta_{jk} = O(m^{-1/2+2\eps}+n^{-1/2+2\eps}) .
\end{align*}
Note that the $O(\,)$ estimates in the last four lines are uniform over
$j,j',k,k'$.

\nicebreak
\subsection{Estimating the main part of the integral}\label{s:complete}

Define $E_2=\exp(L_2)$.  We have shown that the value of the
integrand in
$\Q\cap\M$ is $E_1=E_2\(1+\OO(n^{-1/2})\)$.
Denote the complement of the region $\M$ by $\M^c$.   We can
approximate our integral as follows:
\begin{align}
\int_{\Q \cap \M} E_1
  &= \int_{\Q \cap \M} E_2 + \OO(n^{-1/2})\int_{\Q \cap \M} \abs{E_2} \notag\\[0.4ex]
  &= \int_{\Q \cap \M} E_2 + \OO(n^{-1/2})\int_{\Q} \,\abs{E_2} \notag\\[0.4ex]
  &= \int_{\Q} E_2
        + O(1)\int_{\Q\cap \M^c} \abs{E_2} + \OO(n^{-1/2})\int_{\Q} \,\abs{E_2}.
                  \label{recipe}
\end{align}
It suffices to estimate the value of each integral in \eqref{recipe}.

We first compute the integral of $E_2$ over $\Q$.
We proceed in three stages, starting with integration
with respect to $\mu$.
For the latter, we can use the formula
\[
 \int\limits_{-(mn)^{-1/2+2\eps}}^{(mn)^{-1/2+2\eps}}
   \kern-2mm\exp\(-Amn\mu^2 - i\beta\mu\) d\mu
 = \Bigl(\frac{\pi}{Amn}\Bigr)^{1/2}
   \exp\biggl( -\frac{\beta^2}{4Amn} +O(n^{-1}) \biggr),
\]
provided $\beta=o(A(mn)^{1/2+2\eps})$.  In our case,
$\beta=3A_3(n\mu_2+m\nu_2)$, which is small enough
because of the assumptions $m=o(A^2n^{1+\eps})$
and $n=o(A^2m^{1+\eps})$.
Therefore, integration over $\mu$ contributes
\begin{equation}\label{muint}
\Bigl(\frac{\pi}{Amn}\Bigr)^{1/2}
 \exp\biggl(\frac{-9A_3^2(n\mu_2+m\nu_2)^2}{4Amn}
+O(n^{-1}) \biggr).
\end{equation}

The second step is to
integrate with respect to $\sigmavec$ the integrand
\begin{align}\label{sigbits}
\begin{split}
   \exp\biggl(\! &- An\mu_2
         - \jsum \alpha_{j\c} \sigma_j^2
         - \frac{9A_3^2n}{4Am}\mu_2^2
          - i A_3 n\mu_3
               - 3i A_3 cn \mu_1\mu_2  \\
         &~{}  - i \jsum \beta_{j\c}\sigma_j^3
                 + i\sum_{j,j'=1}^{m-1} g_{jj'}\sigma_j\sigma_{j'}^2
             + i\jsum\ksum
                \(u_{jk}\sigma_j\tau_k^2 + v_{jk}\sigma_j^2\tau_k\) \\
       &~{}+ \Bigl(6A_4-\frac{9A_3^2}{2A}\Bigr)\mu_2\nu_2 + A_4n\mu_4
       + O(n^{-1})\biggr).
\end{split}
\end{align}
This is accomplished by an appeal to Theorem~\ref{MW3},  presented
in the Appendix.  In the terminology of that theorem, we have
$N = m-1$, $\delta(N)=O(n^{-1})$, 
$\eps'=\tfrac32\eps$, $\eps''=\tfrac53\eps$, $\eps'''=3\eps$,
$\bar\eps=6\eps$, and $\hat\eps(N)=\eps+o(1)$ is defined
by $2n^{-1/2+\eps}=N^{-1/2+\hat\eps}$.
Furthermore,
\begin{align*}
\mwA &= \frac{An}{m-1}, &
\mwa_j &=  -\alpha_{j\c}
                     + \Bigl(6A_4-\frac{9A_3^2}{2A}\Bigr)\nu_2
                   + i\ksum v_{jk}\tau_k, &&\\
\mwB_j &= -\frac{iA_3n}{m-1}  - \frac{i}{m-1}\beta_{j\c}, &
\mwC_{jj'} &= -3iA_3cn + i g_{jj'}, \\
\mwE_j &= \frac{A_4n}{m-1}, &
\mwF_{jj'} &= -\frac{9A_3^2n}{4Am}, \\
\mwJ_j &= i\ksum u_{jk}\tau_k^2\,.
\end{align*}
We can take $\Deltait=\tfrac34$, and
calculate that
\begin{gather}
 \frac{3}{4\mwA^2N}\sum_{j=1}^N \mwE_j
    + \frac{1}{4\mwA^2N^2}\sum_{j,j'=1}^N \mwF_{jj'} = \frac{m}{n}
  \biggl(\frac{3A_4}{4A^2}-\frac{9A_3^2}{16A^3}\biggr) + \OO(n^{-1}),
 \notag\displaybreak[0] \\
 \frac{15}{16\mwA^3N}\sum_{j=1}^{N}\mwB_j^2
     + \frac{3}{8\mwA^3N^2}\sum_{j,j'=1}^{N}\mwB_j\mwC_{jj'}
  + \frac{1}{16\mwA^3N^3}\sum_{j,j'\!,j''=1}^N \mwC_{jj'}\mwC_{jj''}
          =  -\frac{3A_3^2m}{8A^3n} + \OO(n^{-1}),
 \notag\displaybreak[0] \\
    \begin{split}
          \frac{1}{2\mwA N}\sum_{j=1}^{N}\mwa_j
     + \frac{1}{4\mwA^2N^2}\sum_{j=1}^{N}\mwa_j^2
      &=  - \frac{1}{2An}\alpha_{\c\c}
      + \frac{1}{4A^2n^2}\jsum (\alpha_{j\c})^2 \\
   &\quad{} + \frac{m}{n}\biggl(\frac{3A_4}{A}-\frac{9A_3^2}{4A^2}\biggr)\nu_2
              + \frac{i}{2An}\ksum v_{\c k}\tau_k+\OO(n^{-1/2}),
    \end{split}\label{taubit} \\
    \mwZ = Z_1 = \exp\biggl( \frac{3A_3^2m}{8A^3n} + \OO(n^{-1}) \biggr)
         = O(1) \exp\biggl(\frac {(1-2\lambda)^2 m}{24 An}\biggr).\notag
\end{gather}

Applying Theorem~\ref{MW3}, we see that  $\Thetait_2 = \OO(n^{-1/2})$,
and so integration with respect to~$\sigmavec$
contributes a $\tau$-free factor
\begin{align}
 \begin{split}\label{contrib1}
 \Bigl(\frac{\pi}{An}\Bigr)^{\!(m-1)/2}
 \exp\biggl(&
  \frac{m}{n} \Bigl(\frac{3A_4}{4A^2}-\frac{15A_3^2}{16A^3}\Bigr)
  - \frac{1}{2A n}\alpha_{\c\c} \\
  &{}+ \frac{1}{4 A^2 n^2}\jsum (\alpha_{j\c})^2
    + \OO(n^{-1/2})+ O(n^{-3/4}Z_1)
    \biggr).
\end{split}
\end{align}
By the conditions of Theorem~\ref{bigtheorem}, $Z_1\le n^{1/5}$,
so $\OO(n^{-1/2})+ O(n^{-3/4}Z_1)=\OO(n^{-1/2})=o(1)$ as required by Theorem~\ref{MW3}.

Finally, we need to integrate over $\tauvec$.  Collecting the
remaining terms from \eqref{L2value}, and the terms
involving $\tauvec$ from
\eqref{muint} and \eqref{taubit}, we have an
integrand equal to
\begin{align*}
&
\exp\biggl(
-Am\nu_2
+\Bigl(\frac{3A_4m}{An}-\frac{9A_3^2m}{4A^2n}\Bigr)\nu_2
-\frac{9A_3^2m}{4An}\nu_2^2
+A_4m\nu_4  -iA_3m\nu_3 -3iA_3dm\nu_2\nu_1\\
& \kern11mm{}
- \ksum \alpha_{\c k}\tau_k^2
- i\ksum\beta_{\c k}\tau_k^3
+\frac{i}{2An}\ksum v_{\c k}\tau_k
 + i\sum_{k,k'=1}^{n-1} h_{kk'}\tau_k\tau_{k'}^2 + \OO(n^{-1/2})
\biggr).
\end{align*}

In the terminology of Theorem~\ref{MW3},
$N = n-1$, $\delta(N)=\OO(n^{-1/2})$,
$\eps'=\tfrac32\eps$, $\eps''=\tfrac53\eps$, $\eps'''=3\eps$,
$\bar\eps=4\eps$, and $\hat\eps(N)=\eps+o(1)$
is defined by $2m^{-1/2+\eps}=N^{-1/2+\hat\eps}$.
Furthermore,
\begin{align*}
\mwA &= \frac{Am}{n-1}, &
\mwa_k &= \frac{3A_4m}{An}-\frac{9A_3^2m}{4A^2n}
                  -\alpha_{\c k}, && \\
\mwB_k &= -\frac{iA_3m}{n-1}  - \frac{i}{n-1}\beta_{\c k} ,&
\mwC_{kk'} &= -3iA_3dm + i h_{kk'}, \\
\mwE_k &= \frac{A_4m}{n-1}, &
\mwF_{kk'} &= -\frac{9A_3^2m}{4An}, \\
\mwJ_k &= \frac{i}{2An}v_{\c k}\,.
\end{align*}
We can take $\Deltait=\frac34$ again and
calculate that
\begin{gather}
 \frac{3}{4\mwA^2N}\sum_{k=1}^N \mwE_k
    + \frac{1}{4\mwA^2N^2}\sum_{k,k'=1}^N \mwF_{kk'} = \frac{n}{m}
  \biggl(\frac{3A_4}{4A^2}-\frac{9A_3^2}{16A^3}\biggr) + \OO(n^{-1}), \notag\\
 \frac{15}{16\mwA^3N}\sum_{k=1}^{N}\mwB_k^2
     + \frac{3}{8\mwA^3N^2}\sum_{k,k'=1}^{N}\mwB_j\mwC_{kk'}
  + \frac{1}{16\mwA^3N^3}\!\sum_{k,k'\!,k''=1}^N \mwC_{kk'}\mwC_{kk''}
          =  -\frac{3A_3^2n}{8A^3m} + \OO(n^{-1}), \notag\\
    \begin{split}
          \frac{1}{2\mwA N}\sum_{k=1}^{N}\mwa_k
     + \frac{1}{4\mwA^2N^2}\sum_{k=1}^{N}\mwa_k^2
      &=  - \frac{1}{2Am}\alpha_{\c\c}
      + \frac{1}{4A^2m^2}\ksum (\alpha_{\c k})^2 \\
         &\quad{} - \frac{9A_3^2}{8A^3} + \frac{3A_4}{2A^2} +\OO(n^{-1/2}), \\
    \end{split}\notag \\
    \mwZ = Z_2 = \exp\biggl( \frac{3A_3^2n}{8A^3m} + \OO(n^{-1}) \biggr)
         = O(1) \exp\biggl(\frac {(1-2\lambda)^2 n}{24 Am}\biggr).\notag
\end{gather}

We again find that $\Thetait_2=\OO(n^{-1/2})$.  Including the contributions
from \eqref{muint} and \eqref{contrib1}, we obtain
\begin{equation}\label{QE2a}
\begin{split}
 \int_\Q E_2
  &= \Bigl(\frac{\pi}{Amn}\Bigr)^{1/2}
        \Bigl(\frac{\pi}{An}\Bigr)^{\!(m-1)/2}
        \Bigl(\frac{\pi}{Am}\Bigr)^{\!(n-1)/2} \\
  &\qquad{}\times\exp\biggl(  - \frac{9A_3^2}{8A^3} + \frac{3A_4}{2A^2}
         +\Bigl(\frac mn+\frac nm\Bigr)
         \Bigl( \frac{3A_4}{4A^2}-\frac{15A_3^2}{16A^3}\Bigr) \\
   &\kern23mm{} - \Bigl(\frac{1}{2Am}+\frac{1}{2An}\Bigr)\alpha_{\c\c}
      + \frac{1}{4A^2m^2}\ksum (\alpha_{\c k})^2 \\
   &\kern23mm{}
      + \frac{1}{4A^2n^2}\jsum (\alpha_{j\c})^2
      + \OO(n^{-1/2}) Z_2 \biggr).
\end{split}
\end{equation}

Using \eqref{rad10} and the conditions of
Theorem~\ref{bigtheorem}, we calculate that
\begin{align*}
\alpha_{\c\c} &= -\frac12
                      \Bigl(\frac{R}{n}+\frac{C}{m}\Bigr)+\OO(n^{1/2}), \\
\jsum (\alpha_{j\c})^2 &= \dfrac14 (1-2\lambda)^2R+\OO(n^{3/2}), \\
\ksum (\alpha_{\c k})^2 &= \dfrac14 (1-2\lambda)^2 C+\OO(n^{3/2}), \\
\OO(n^{-1/2}) Z_2 &= \OO(n^{-1/2}) n^{2a/5} = O(n^{-b}).
\end{align*}
Substituting these values into \eqref{QE2a} together with the
actual values of $A, A_3, A_4$, we conclude that
\begin{equation}\label{QE2}
\begin{split}
 \int_\Q E_2
  &= \Bigl(\frac{\pi}{Amn}\Bigr)^{1/2}
        \Bigl(\frac{\pi}{An}\Bigr)^{\!(m-1)/2}
        \Bigl(\frac{\pi}{Am}\Bigr)^{\!(n-1)/2} \\
  &\qquad{}\times\exp\biggl(  - \frac12
         - \frac{1-2A}{24A}\Bigl(\frac mn+\frac nm\Bigr)
   + \frac{1}{4A}\Bigl(\frac{1}{m}+\frac{1}{n}\Bigr)
               \Bigl(\frac{R}{n}+\frac{C}{m}\Bigr) \\
      &\kern23mm{}
      + \frac{1-8A}{16A^2}\Bigl(\frac{R}{n^2}+\frac{C}{m^2}\Bigr)
      + O(n^{-b}) \biggr).
\end{split}
\end{equation}

\medskip
We next infer a estimate of $\int_\Q\, \abs{E_2}$.  The calculation that
lead to \eqref{QE2a} remains valid if we set all the values
$A_3$, $\beta_{jk}$, $g_{jj'}$, $h_{kk'}$, $u_{jk}$ and $v_{jk}$
to zero, which is the same as replacing $L_2$ by its real part.
Since $\abs{E_2} = \exp\(\Re(L_2)\)$, this gives
\begin{align}
  \int_\Q\, \abs{E_2}
    &= \exp\biggl(\frac{9A_3^2}{8A^3} 
         + \frac{15A_3^2}{16A^3}\Bigl(\frac mn+\frac nm\Bigr)
         + o(1)\biggr)\int_Q E_2 \notag\\[0.3ex]
    &= \exp\biggl(\frac{(1-2\lambda)^2}{8A}\Bigl( 1 +
         \frac{5n}{6m} + \frac{5m}{6n}\Bigr)+ o(1)\biggr)\int_Q E_2 \notag\\
     &= O(n^a)\int_Q E_2 \label{absE2}
\end{align}
under the assumptions of Theorem~\ref{bigtheorem}.
The third term of \eqref{recipe} can now be identified:
\begin{equation}\label{OEbit}
   \OO(n^{-1/2})  \int_\Q\, \abs{E_2} = \OO(n^{-1/2}) n^a \int_Q E_2
        = O(n^{-b})  \int_Q E_2,
\end{equation}
where, as always, we suppose that $\eps$ is sufficiently small.

Finally, we consider the second term of \eqref{recipe}, namely
\[\int_{\Q\cap \M^c} \abs{E_2},\]
which we will bound as a fraction of $\int_\Q\, \abs{E_2}$ using
a statistical technique.  The following is a well-known result of
Hoeffding~\cite{hoeffding}.

\begin{lemma}\label{hoeffding}
  Let $X_1,X_2,\ldots,X_N$ be independent random
variables such that $\expect X_i=0$ and $\abs{X_i}\le M$ for all~$i$.
Then, for any $t\ge 0$,
\[\Prob\Bigl(\,\sum_{i=1}^N X_i \ge t\Bigr)
    \le \exp\biggl(-\frac{t^2}{2NM^2}\biggr).\]
\end{lemma}

Now consider $\abs{E_2} = \exp\(\Re(L_2)\)$.
Write $\M=\M_1\cap\M_2$, where
$\M_1=\{\, \abs{\mu_1}\le m^{1/2} n^{-1/2+3\eps}\,\}$
and $\M_2=\{\, \abs{\nu_1}\le n^{1/2} m^{-1/2+3\eps}\}$.
For fixed values
of $\mu$ and $\sigmavec$, $\Re(L_2)$ separates over
$\tau_1,\tau_2,\ldots,\tau_{n-1}$ and therefore, apart from
normalization,  it is the joint density of independent
random variables $X_1,X_2,\ldots,X_{n-1}$ which satisfy
$\expect X_k=0$ (by symmetry) and $\abs{X_k}\le 2m^{-1/2+\eps}$
(by the definition of $\Q$).  By Lemma~\ref{hoeffding}, the
fraction of the integral over $\tauvec$ (for fixed $\mu,\sigmavec$)
that has
$\nu_1\ge n^{1/2}m^{-1/2+3\eps}$ is at most
$\exp(-m^{4\eps}/2)$.  By symmetry, the same bound holds
for $\nu_1\le -n^{1/2}m^{-1/2+3\eps}$.  Since these bounds
are independent of $\mu$ and $\sigmavec$, we have
\[\int_{\Q\cap\M_2^c} \,\abs{E_2}
   \le 2\exp(-m^{4\eps}/2)\int_\Q\,\abs{E_2}.\]
By the same argument,
\[\int_{\Q\cap\M_1^c} \,\abs{E_2}
   \le 2\exp(-n^{4\eps}/2)\int_\Q\,\abs{E_2}.\]
Therefore we have in total that
\begin{align}\label{corners}
\int_{\Q\cap\M^c} \,\abs{E_2}
   &\le 2\(\exp(-m^{4\eps}/2)+\exp(-n^{4\eps}/2)\) \int_\Q\,\abs{E_2}
      \notag \\
   &\le  O(n^{-b})  \int_Q E_2,
\end{align}
as for~\eqref{OEbit}.  Applying \eqref{recipe} with \eqref{QE2},
\eqref{OEbit} and~\eqref{corners}, we find that $\int_{\Q\cap\M} E_1$
is given by~\eqref{QE2}.
Multiplying by the Jacobians of the transformations~$T_2$ and~$T_3$,
we find that Theorem~\ref{Jintegral} is proved for $\S'$
given by~\eqref{RSprime}.

\section{Bounding the remainder of the integral}\label{s:boxing}

In the previous section, we estimated the value of the integral
$I_{\R'}(\svec,\tvec)$, which is the same as $I(\svec,\tvec)$
except that it is restricted to a certain region $\R'\supseteq\R$
(see (\ref{FDef}--\ref{bigR})) In this section, we
extend this to an estimate of $I(\svec,\tvec)$ by showing that
the remainder of the region of integration contributes
negligibly.

Precisely, we show the following.
\begin{thm}\label{boxing}
  Let $F(\thetavec,\phivec)$ be the integrand of $I(\svec,\tvec)$
  as defined in~\eqref{FDef}.
  Then, under the conditions of Theorem~\ref{bigtheorem},
  \[
   \int_{\R^c} \abs{F(\thetavec,\phivec)}\,d\thetavec d\phivec
    = O(n^{-1}) \int_{\R'} F(\thetavec,\phivec)\,d\thetavec d\phivec.
  \]
\end{thm}

For $1\leq j\leq m$, $1\leq k\leq n$, let
$A_{jk} = A + \alpha_{jk}=\tfrac12\lambda_{jk}(1-\lambda_{jk})$
(recall \eqref{AlBetGamDef}),
and define $\Amin=\min_{jk} A_{jk} = A+\OO(n^{-1/2})$.
We begin with two technical lemmas whose proofs are omitted.
\begin{lemma}\label{fbnd}
\[\abs{F(\thetavec,\phivec)} = \prod_{j=1}^m\prod_{k=1}^n f_{jk}(\t_j + \p_k),\]
where
\[f_{jk}(z) =\ssqrt{1-4A_{jk}(1-\cos z)}\,.\]
Moreover, for all real $z$,
\[0\le f_{jk}(z) \le \exp\( -A_{jk} z^2 + \dfrac1{12} A_{jk} z^4\).
    \quad\qedsymbol\]
\end{lemma}

\begin{lemma}\label{ibnd}
For all $c>0$,
\[\int_{-8\pi/75}^{8\pi/75} \exp\( c (-x^2 + \dfrac73 x^4)\)\,dx
 \le \sqrt{\pi/c}\, \exp(3/c).   \quad\qedsymbol\]
\end{lemma}

\begin{proof}[Proof of Theorem \ref{boxing}]
Our approach will be to bound $\int \,\abs{F(\thetavec,\phivec)}$ over
a variety of regions whose union covers $\R^c$.
To make the comparison of these bounds with
$\int_{\R'} F(\thetavec,\phivec)$ easier, we note that
\begin{equation}\label{I0I1}
   \int_{\R'} F(\thetavec,\phivec)\,d\thetavec d\phivec
     = \exp\(O(m^\eps+n^\eps)\) I_0 = \exp\(O(m^{3\eps}+n^{3\eps})\) I_1,
\end{equation}
where
\begin{align*}
   I_0 &= \Bigl(\frac{\pi}{A_{\b\b}}\Bigr)^{1/2}\,
   \prod_{j=1}^m   \Bigl(\frac{\pi}{A_{j\b} }\Bigr)^{\!1/2}\,
          \prod_{k=1}^n   \Bigl(\frac{\pi}{A_{\b k} }\Bigr)^{\!1/2}
          \!, \\
   I_1 &= \Bigl(\frac{\pi}{An}\Bigr)^{\!m/2} \Bigl(\frac{\pi}{Am}\Bigr)^{\!n/2}\!.
\end{align*}
To see this, expand 
\[
 A_{j\b} = An+\alpha_{j\b}=An \exp\biggl( \frac{\alpha_{j\b}}{An}
   - \frac{\alpha_{j\b}^2}{2A^2n^2} +\cdots\,\biggr),\]
and similarly for $A_{\b k}$, and compare the result
to Theorem~\ref{Jintegral} using the assumptions of
Theorem~\ref{bigtheorem}.  It may help to recall the calculation
following~\eqref{QE2a}.

\smallskip

Take $\kappa=\pi/300$ and define $x_0,x_1,\ldots,x_{299}$ by
$x_\ell = 2\ell\kappa$.
For any $\ell$,
let $\S_1(\ell)$ be the set of
$(\thetavec,\phivec)$ such that $\t_j\in [x_\ell-\kappa,x_\ell+\kappa]$ for
at least $\kappa m/\pi$ values of~$j$ and
$\p_k\notin [-x_\ell-2\kappa,-x_\ell+2\kappa]$ for at
least $n^\eps$ values of~$k$.
For $(\thetavec,\phivec)\in\S_1(\ell)$, $\theta_j+\phi_k\notin[-\kappa,\kappa]$
for at least $\kappa mn^\eps/\pi$ pairs $(j,k)$ so,
by Lemma~\ref{fbnd},
$\abs{F(\thetavec,\phivec)}\le \exp(-c_1\Amin mn^\eps)$
for some $c_1>0$ which is independent of~$\ell$.

Next  define $\S_2(\ell)$ to be the set of
$(\thetavec,\phivec)$ such that
$\t_j\in [x_\ell-\kappa,x_\ell+\kappa]$ for
at least $\kappa m/\pi$ values of~$j$,
$\p_k\in[-x_\ell-2\kappa,-x_\ell+2\kappa]$ for at least $n-n^\eps$ values of~$k$
and
$\t_j\notin [x_\ell-3\kappa,x_\ell+3\kappa]$ for at least $m^\eps$ values of~$j$.
By the same argument with the roles of $\thetavec$
and $\phivec$ reversed,
$\abs{F(\thetavec,\phivec)}\le \exp(-c_2\Amin m^\eps n)$
for some $c_2>0$ independent of~$\ell$ when
$(\thetavec,\phivec)\in\S_2(\ell)$.

Now define $\R_1(\ell)$ to be the set of
pairs $(\thetavec,\phivec)$ such that
$\t_j\in [x_\ell-3\kappa,x_\ell+3\kappa]$ for at least $m-m^\eps$
values of~$j$, and $\p_k\in [-x_\ell-3\kappa,-x_\ell+3\kappa]$ for
at least $n-n^\eps$ values of~$k$.
By the pigeonhole principle, for any $\thetavec$ there is some~$\ell$
such that $[x_\ell-\kappa,x_\ell+\kappa]$ contains at least $\kappa m/\pi$
values of $\t_j$.  Therefore,
\[
    \Bigl(\,\bigcup_{\ell=0}^{299} \, \R_1(\ell)\!\Bigr)^c 
        \subseteq \bigcup_{\ell=0}^{299}  \,\(\S_1(\ell)\cup\S_2(\ell)\).
\]
Since the total volume of $\(\,\bigcup_\ell\R_1(\ell)\)^c$ is at most $(2m)^{m+n}$,
we find that for some $c_3>0$,
\begin{align}
 \int_{(\bigcup_\ell \R_1(\ell))^c} \,&
    \abs{F(\thetavec,\phivec)}\,d\thetavec d\phivec \notag\\
  & \le
     (2\pi)^{m+n} \( \exp( -c_3 \Amin m n^\eps)
         +\exp( -c_3 \Amin m^\eps n)\) \notag\\
  & \le e^{-n} I_1.\label{R1C}
\end{align}

We are left with $(\thetavec,\phivec)\in\bigcup_\ell\R_1(\ell)$.
If we subtract $x_\ell$ from each $\t_j$ and add $x_\ell$ to each $\p_k$ the
integrand $F(\thetavec,\phivec)$ is unchanged, so we can assume
for convenience that $\ell=0$ and that
$(\thetavec,\phivec)\in\R_1=\R_1(0)$.
The bounds we obtain on parts of the integral we seek to reject will
be at least $1/300$ of the total and thus be of the right order of
magnitude.  We will not mention this point again.

For a given $\thetavec$, partition $\{1,2,\ldots, m\}$ into sets
$J_0=J_0(\thetavec)$, $J_1 = J_1(\thetavec)$ and $J_2 = J_2(\thetavec)$,
containing the indices $j$ such that
$\abs{\t_j} \leq 3\kappa$, $3\kappa  < \abs{\t_j}\leq 15\kappa$
and $\abs{\t_j} > 15\kappa$, respectively.
Similarly partition $\{1,2,\ldots,n\}$ into
$K_0=K_0(\phivec)$, $K_1=K_1(\phivec)$ and $K_2=K_2(\phivec)$.
The value of $\abs{F(\thetavec,\phivec)}$ can now be bounded using
\begin{align*}
f_{jk}(\t_j+\p_k) & \\
 &\kern-14mm{}\le \begin{cases}
    \,\exp\(-\Amin(\t_j+\p_k)^2+\tfrac1{12}\Amin(\t_j+\p_k)^4\)
        & \textrm{if }(j,k) \in (J_0\cup J_1)\times(K_0\cup K_1), \\[0.8ex]
    \,\ssqrt{\vrule width0pt depth0.3ex
          1-4\Amin(1-\cos(12\kappa))} \le e^{-\Amin/64}
        & \textrm{if }(j,k) \in (J_0\times K_2)\cup (J_2\times K_0), \\[0.3ex]
    \,1   & \textrm{otherwise}.
 \end{cases}
\end{align*}
Let $I_2(m_2,n_2)$ be the contribution to
$\int_{\R_1} \abs{F(\thetavec,\phivec)}$
of those $(\thetavec,\phivec)$ with
$\card{J_2}=m_2$ and $\card{K_2}=n_2$.
Recall that $\card{J_0}> m-m^\eps$ and $\card{K_0}> n-n^\eps$.
We have
\begin{align}\label{I2bnd}
\begin{split}
  I_2(m_2,n_2) 
 &\le \binom{m}{m_2}\binom{n}{n_2}(2\pi)^{m_2+n_2} \\
       &\kern7mm{}\times \exp\( -\dfrac1{64} \Amin(n-n^\eps)m_2
               -\dfrac1{64} \Amin(m-m^\eps)n_2\) I'_2(m_2,n_2),
\end{split}
\end{align}
where
\[
  I'_2(m_2,n_2) =
   \int_{-15\kappa}^{15\kappa}\!\!\cdots \int_{-15\kappa}^{15\kappa}
    \!\!\!\exp\Bigl(- \Amin\sumpp_{jk}(\t_j+\p_k)^2 +
                \dfrac1{12}\Amin\sumpp_{jk}(\t_j+\p_k)^4\Bigr)
      \, d\thetavec' d\phivec',
\]
and the primes denote restriction to $j\in J_0\cup J_1$
and $k\in K_0\cup K_1$.
Write $m'=m-m_2$ and $n'=n-n_2$ and define
$\avtheta'=(m')^{-1}\sumpp_j\t_j$,
$\breve\theta_j=\t_j-\avtheta'$ for $j\in J_0\cup J_1$,
$\avphi'=(n')^{-1}\sumpp_k\p_k$,
$\breve\phi_k=\p_k-\avphi'$ for $k\in K_0\cup K_1$,
$\nu'=\avphi'+\avtheta'$ and
$\delta'=\avtheta'-\avphi'$.
Change variables from $(\thetavec',\phivec')$ to
$\{\breve\theta_j \mathrel| j\in J_3\} 
\cup \{\breve\phi_k\mathrel| k\in K_3\}
\cup \{\nu',\delta'\}$, where $J_3$ is some subset of $m'-1$ elements
of $J_0\cup J_1$ and $K_3$ is some subset of $n'-1$ elements
of $K_0\cup K_1$.  {}From Section~\ref{s:integral} we know that the
Jacobian of this transformation is $m'n'/2$.
The integrand of $I'_2$ can now be bounded using
\[ \sumpp_{jk}(\t_j+\p_k)^2 =
   n'\sumpp_j\breve\theta_j^2 + m'\sumpp_k\breve\phi_k^2
   + m'n'\nu'{}^2
\]
and
\[ \sumpp_{jk}(\t_j+\p_k)^4 \le 
     27n'\sumpp_j\breve\theta_j^4 + 27m'\sumpp_k\breve\phi_k^4
     + 27m'n'\nu'{}^4.
\]
The latter follows from the inequality
$(x+y+z)^4 \le 27(x^4+y^4+z^4)$ valid for all $x,y,z$.
Therefore,
\begin{align*}
I'_2(m_2,n_2) \le \frac{O(1)}{m'n'} \int_{-30\kappa}^{30\kappa}
    \int_{-30\kappa}^{30\kappa}\!\!\cdots \int_{-30\kappa}^{30\kappa}
     \exp\Bigl(& \Amin n' \sumpp_jg(\breve\theta_j) 
    + \Amin m'\sumpp_k g(\breve\phi_k) \\
    &{}+ \Amin m'n' g(\nu') \Bigr)\,
    d\breve\theta_{j\in J_3}d\breve\phi_{k\in K_3}\, d\nu',
\end{align*}
where $g(z) = -z^2+\tfrac94 z^4$.  Since $g(z)\le 0$ for
$\abs{z}\le 30\kappa$, and we only need
an upper bound, we can restrict the summations in the
integrand to $j\in J_3$ and
$k\in K_3$.  The integral now separates into $m'+n'-1$
one-dimensional integrals
and Lemma~\ref{ibnd} (by monotonicity) gives that
\[
I'_2(m_2,n_2) = O(1)
     \frac{\pi^{(m'+n')/2}}{ \Amin^{(m'+n'-1)/2} (m')^{n'/2-1} (n')^{m'/2-1}}
     \exp\( O(m'/(\Amin n') + n'/(\Amin m'))\).
\]
Applying \eqref{I0I1} and \eqref{I2bnd}, we find that
\begin{equation}\label{I2mn}
 \mathop{\sum_{m_2=0}^{m^\eps} \;
                \sum_{n_2=0}^{n^\eps}}\displaylimits_{m_2+n_2\ge 1}
 I_2(m_2,n_2) = O\(e^{-c_4A m}+e^{-c_4A n}\) I_1
\end{equation}
for some $c_4>0$.

\medskip

We have now bounded contributions to the integral of
$\abs{F(\thetavec,\phivec)}$ from everywhere outside
the region
\[
\X = \bigl\{\, (\thetavec,\phivec) \bigm|
  \abs{\t_j},\abs{\p_k} \le15\kappa~\text{for}~1\le j\le m, 1\le k\le n\, \bigr\}.
\]
By Lemma~\ref{fbnd}, we have for $(\thetavec,\phivec)\in\Chat^{m+n}$
(which includes~$\X$) that
\[\abs{F(\thetavec,\phivec)} \le\exp\Bigl(
 - \sum_{j=1}^m\sum_{k=1}^n A_{jk} (\that_j+\phat_k+\nu)^2
  + \dfrac1{12}\sum_{j=1}^m\sum_{k=1}^n A_{jk} (\that_j+\phat_k+\nu)^4\Bigr),
\]
where $\that_j = \theta_j - \avtheta$, $\phat_k = \phi_k - \avphi$
and $\nu=\avtheta+\avphi$.  As before, the integrand is independent
of $\delta=\avtheta-\avphi$ and our notation will tend to ignore $\delta$
for that reason; for our bounds it will suffice to remember that
$\delta$ has a bounded range.

We proceed by exactly diagonalizing the $(m+n+1)$-dimensional
quadratic form.  Since $\sum_{j=1}^m\that_j = \sum_{k=1}^n\phat_k=0$,
we have
\begin{align*}
 \sum_{j=1}^{m}\sum_{k=1}^{n}
   A_{jk} (\that_j+\phat_k+\nu)^2 &= \sum_{j=1}^{m} A_{j\b}\that_j^2
    + \sum_{k=1}^n A_{\b k}\phat_k^2 + A_{\b\b}\nu^2\\
     &\quad{}+ 2\sum_{j=1}^m\sum_{k=1}^n \alpha_{jk}\that_j\phat_k
           + 2\nu\sum_{j=1}^m \alpha_{j\b}\that_j + 2\nu\sum_{k=1}^n
                            \alpha_{\b k}\phat_k.
\end{align*}
This is almost diagonal, because $\alpha_{jk} = \OO(n^{-1/2})$, and we can
correct it with the slight additional transformation $(I+X^{-1}Y)^{-1/2}$
described by Lemma~\ref{diagonalize}, where $X$ is a diagonal
matrix with diagonal entries $A_{j\b}$, $A_{\b k}$ and~$A_{\b\b}$.
The matrix $Y$ has zero diagonal and other entries of
magnitude $\OO(n^{-1/2})$ apart from the row and column indexed by $\nu$, which
have entries of magnitude~$\OO(n^{1/2})$.
By the same argument as used in Section~\ref{s:diagonalize},
all eigenvalues of $X^{-1}Y$ have  magnitude~$\OO(n^{-1/2})$, so the
transformation is well-defined.
The new variables
$\{\hat\vartheta_j\}$, $\{\hat\varphi_k\}$ and $\varnu$
are related to the old by
\[
(\that_1,\ldots,\that_m,\phat_1,\ldots,\phat_n,\nu)^T
 =(I+X^{-1}Y)^{-1/2}(\hat\vartheta_1,\ldots,\hat\vartheta_m,
 \hat\varphi_1,\ldots,\hat\varphi_n,\varnu)^T.
\]
We will keep the variable~$\delta$ as a variable of
integration but, as noted before, our notation will generally ignore~it.

More explicitly, for some $d_1,\ldots,d_m,d'_1,\ldots,d'_n=\OO(n^{-3/2})$,
we have uniformly over\\ $j=1,\ldots,m$, $k=1,\ldots,n$ that
\begin{align}\label{varrels}
\begin{split}
  \that_j &= \hat\vartheta_j + \sum_{q=1}^m \OO(n^{-2})\hat\vartheta_q
                         + \sum_{k=1}^n \OO(n^{-3/2})\hat\varphi_k + \OO(n^{-1/2})\varnu, \\
  \phat_k &= \hat\varphi_k + \sum_{j=1}^m \OO(n^{-3/2})\hat\vartheta_j
                         + \sum_{q=1}^n \OO(n^{-2})\hat\varphi_q + \OO(n^{-1/2})\varnu, \\
  \nu &= \varnu + \sum_{j=1}^m d_j\hat\vartheta_j
                         + \sum_{k=1}^n d'_k \hat\varphi_k + \OO(n^{-1})\varnu.
\end{split}
\end{align}
Note that the expressions $O(\,)$ in \eqref{varrels} represent values
that depend on $m,n,\svec,\tvec$ but not on $\{\hat\vartheta_j\},
\{\hat\varphi_k\},\varnu$.

The region of integration $\X$ is
$(m{+}n)$-dimensional.  In place of the variables
$(\thetavec,\phivec)$ we can use $(\thetahatvec,\phihatvec,\nu,\delta)$
by applying the identities $\that_m=-\sum_{j=1}^{m-1}\that_j$ and
$\phat_n=-\sum_{k=1}^{n-1}\phat_k$.
(Recall that $\thetahatvec$ and $\phihatvec$ don't include
$\that_m$ and $\phat_n$.)
The additional
transformation \eqref{varrels} maps the two just-mentioned identities
into identities that define $\hat\vartheta_m$ and $\hat\varphi_n$
in terms of $(\varthetahatvec,\varphihatvec,\varnu)$, where
$\varthetahatvec=(\hat\vartheta_1,\ldots,\hat\vartheta_{m-1})$
and $\varphihatvec=(\hat\varphi_1,\ldots,\hat\varphi_{n-1})$.
These have the form
\begin{equation}\label{mnvars}
\begin{split}
 \hat\vartheta_m &= -\sum_{j=1}^{m-1}\(1+\OO(n^{-1})\)\hat\vartheta_j
      + \sum_{k=1}^{n-1}\OO(n^{-1/2})\hat\varphi_k + \OO(n^{1/2})\varnu,\\
 \hat\varphi_n &= \sum_{j=1}^{m-1}\OO(n^{-1/2})\hat\vartheta_j 
       - \sum_{k=1}^{n-1}\(1+\OO(n^{-1})\)\hat\varphi_k + \OO(n^{1/2})\varnu.
\end{split}
\end{equation}
Therefore, we can now integrate over
$(\varthetahatvec,\varphihatvec,\varnu,\delta)$.
The Jacobian of the transformation from $(\thetavec,\phivec)$
to $(\thetahatvec,\phihatvec,\nu,\delta)$
is $mn/2$, as in Section~\ref{s:integral}.
The Jacobian of the transformation
$T_4(\varthetahatvec,\varphihatvec,\varnu)
 = (\thetahatvec,\phihatvec,\nu)$ defined
by~\eqref{varrels} can be seen to be $1+\OO(n^{-1/2})$ by
Lemma~\ref{l:detexp}, using the fact that the $\infty$-norm
of the matrix of partial derivatives is~$\OO(n^{-1/2})$.  This matrix
has order $m+n-1$ and can be obtained by substituting
\eqref{mnvars} into~\eqref{varrels}.

The transformation $T_4$ changes the region of integration only by
a factor $1+\OO(n^{-1/2})$ in each direction, since the inverse of \eqref{varrels}
has exactly the same form except that the constants $\{ d_j\}, \{d'_k\}$,
while still of magnitude $\OO(n^{-3/2})$, may be different.
Therefore, the image of
region $\X$ lies inside the region
\[
  \Y = \bigl\{\, (\varthetahatvec,\varphihatvec,\varnu) \bigm|
  \abs{\hat\vartheta_j},\abs{\hat\varphi_k}\le 31\kappa~(1\le j\le m, 1\le k\le n),
   \,\abs\varnu\le 31\kappa \, \bigr\}.
\]

We next bound the value of the integrand in $\Y$.   By repeated application
of the inequality $xy\le \tfrac12 x^2+\tfrac12 y^2$, we find that
\[
\dfrac1{12} \sum_{j=1}^m\sum_{k=1}^n A_{jk} (\that_j+\phat_k+\nu)^4
 \le \dfrac73  \Bigl(\, \sum_{j=1}^m A_{j\b}\hat\vartheta_j^4
    + \sum_{k=1}^n A_{\b k}\hat\varphi_k^4
    + A_{\b\b}\varnu^4\Bigr),
\]
where we have chosen $\tfrac73$ as a convenient value greater
than $\tfrac94$.  Now define $h(z) = -z^2+\tfrac73 z^4$.  Then,
for $(\varthetahatvec,\varphihatvec,\varnu)\in\Y$, 
\begin{align}
  \abs{F(\thetavec,\phivec)}
  &\le \exp\Bigl(\, \sum_{j=1}^m A_{j\b} h(\hat\vartheta_j)
    + \sum_{k=1}^n A_{\b k} h(\hat\varphi_k) + A_{\b\b} h(\varnu) \Bigr)
    \notag\\
    &\le \exp\Bigl(\, \sum_{j=1}^{m-1} A_{j\b} h(\hat\vartheta_j)
    + \sum_{k=1}^{n-1} A_{\b k} h(\hat\varphi_k) + A_{\b\b} h(\varnu) \Bigr)
    \label{I199}\\[-0.5ex]
    &= \exp\(A_{\b\b}h(\varnu)\)
           \prod_{j=1}^{m-1} \exp\(A_{j\b}h(\hat\vartheta_j)\)
          \prod_{k=1}^{n-1} \exp\(A_{\b k}h(\hat\varphi_k)\),\label{I200}
\end{align}
where the second line holds
because $h(z)\le 0$ for $\abs{z}\le 31\kappa$.

Define
 \begin{align*}
   \W_0 &= \bigl\{ (\varthetahatvec,\varphihatvec,\varnu)\in\Y \bigm|
           \abs{\hat\vartheta_j}\le \tfrac12 n^{-1/2+\eps} \; (1\le j\le m-1),\\
           &\kern35mm\abs{\hat\varphi_k}\le \tfrac12 m^{-1/2+\eps} \; (1\le k\le n-1),\\
           &\kern37.5mm\abs{\varnu}\le \tfrac12 (mn)^{-1/2+2\eps}  
          \bigr\}, \\[0.4ex]
    \W_1&=\Y - \W_0, \\
    \W_2 &= \Bigl\{ \, (\varthetahatvec,\varphihatvec,\varnu)\in\Y \Bigm|\,
                      \Bigl| \, \sum_{j=1}^{m-1} d_j\hat\vartheta_j
                              + \sum_{k=1}^{n-1} d'_k\hat\varphi_k\Bigr| \le n^{-5/4} \,\Bigr\}.
\end{align*}
Also define similar regions $\W'_0,\W'_1,\W'_2$ by omitting the variables
$\hat{\vartheta}_1,\hat{\varphi}_1$ instead of
$\hat{\vartheta}_m,\hat{\varphi}_n$ starting at~\eqref{I199}.
Using \eqref{varrels}, we see that $T_4$, and the corresponding
transformation that omits $\hat{\vartheta}_1$ and $\hat{\varphi}_1$,
map $\R$ to a superset of
$\W_0\cap\W_2\cap\W'_1\cap\W'_2$.  Therefore,
$\X-\R$ is mapped to a subset of
$\W_1\cup(\W_0-\W_2)\cup\W'_1\cup(\W'_0-\W'_2)$ and it will
suffice to find a tight bound on the integral in each of the four
latter regions.
 
Denoting the right side of~\eqref{I200} by
$F_0(\varthetahatvec,\varphihatvec,\varnu)$,
Lemma~\ref{ibnd} gives
\begin{equation}\label{I0F0}
\int_\Y
    F_0(\varthetahatvec,\varphihatvec,\varnu)
    \, d\varthetahatvec d\varphihatvec d\varnu
    = \exp\(O(m^\eps+n^\eps)\) I_0.
\end{equation}
Also note that
\begin{equation}\label{1Dtail}
 \int_{z_0}^{31\kappa} \exp( c\, h(z)) = O(1)\exp(c\, h(z_0))
\end{equation}
for $c,z_0>0$ and $z_0=o(1)$, since
$h(z)\le h(z_0)$ for $z_0\le z\le 31\kappa$.
By applying \eqref{1Dtail} to each of the factors of~\eqref{I200} in turn,
\begin{equation}\label{W1bnd}
  \int_{\W_1}  F_0(\varthetahatvec,\varphihatvec,\varnu)
     \, d\varthetahatvec d\varphihatvec d\varnu = 
    O\(e^{-c_6A m^{2\eps} } + e^{-c_6A n^{2\eps} }  \) I_0
\end{equation}
for some $c_6>0$ and so, by \eqref{I0F0} and \eqref{W1bnd},
\[
  \int_{\W_0}  F_0(\varthetahatvec,\varphihatvec,\varnu)
     \, d\varthetahatvec d\varphihatvec d\varnu = 
    \exp\(O(m^\eps+n^\eps)\) I_0.
\]
Applying Lemma~\ref{hoeffding} twice, once to the variables
$d_1\hat\vartheta_1,\ldots,d_{m-1}\hat\vartheta_{m-1},
  d'_1\hat\varphi_1,\ldots,d'_{n-1}\hat\varphi_{n-1}$
and once to their negatives,
using $M=\OO(n^{-2})$, $N=m+n-2$ and $t=n^{-5/4}$,
we find that
\begin{align}
  \int_{\W_0-\W_2}  F_0(\varthetahatvec,\varphihatvec,\varnu)
     \, d\varthetahatvec d\varphihatvec d\varnu
     &= O\( e^{-n^{1/4}}\) \int_{\W_0}  F_0(\varthetahatvec,\varphihatvec,\varnu)
     \, d\varthetahatvec d\varphihatvec d\varnu\notag\\
     &= O\( e^{-n^{1/5}}\) I_0.\label{W0W2}
\end{align}

Finally, parallel computations give the same bounds on the integrals
over $\W'_1$ and $\W'_0-\W'_2$.

We have now bounded
$\int \abs{F(\thetavec,\phivec)}$ in regions that together cover the
complement of~$\R$.  Collecting these bounds from
\eqref{R1C}, \eqref{I2mn}, \eqref{W1bnd},
\eqref{W0W2}, and the above-mentioned analogues of
\eqref{W1bnd} and \eqref{W0W2}, we conclude that
\[
 \int_{\R^c} \abs{F(\thetavec,\phivec)}\,d\thetavec d\phivec
=O\(e^{-c_7A m^{2\eps} } + e^{-c_7A n^{2\eps} }  \) I_0
\]
for some $c_7>0$, which implies the theorem by~\eqref{I0I1}.
\end{proof}

\nicebreak
\section*{Appendix: Estimating an integral}\label{s:appendix}

In this appendix we estimate the value of a certain multi-dimensional integral.
A similar integral appeared in~\cite{MW} and variations of it appeared in
~\cite{Mtourn,euler,MWtournament}.
However, none of the previously published variations
meet our present
requirements entirely.  We will meet them here, and also
introduce a new method of proof that gives a better error term.

It is intended that this appendix be notationally independent of
the rest of the paper.  We have used new symbols where possible,
but even in the few remaining exceptions, assumptions about the
values of variables stated earlier do not apply here.

\begin{thm}\label{MW3}
Let $\eps', \eps'', \eps''', \bar\eps,\Deltait$ be constants
such that $0<\eps'<\eps''<\eps'''$, $\bar\eps\ge 0$,
and $0<\Deltait<1$.  The following is true if $\eps'''$
and $\bar\eps$ are sufficiently small.\endgraf
Let $\mwA=\mwA(N)$ be a real-valued function such that
$\mwA(N)=\Omega(N^{-\eps'})$.
Let $\mwa_j=\mwa_j(N)$, $\mwB_j=\mwB_j(N)$,
$\mwC_{jk}=\mwC_{jk}(N)$, $\mwE_j=\mwE_j(N)$,
$\mwF_{jk}=\mwF_{jk}(N)$ and $\mwJ_j=\mwJ_j(N)$ be
complex-valued functions $(1 \le j,k \le N)$
such that $\mwB_j,\mwC_{jk},
\mwE_j,\mwF_{jk}=O(N^{\bar\eps})$,
$\mwa_j=O(N^{1/2+\bar\eps})$, and $\mwJ_j=O(N^{-1/2+\bar\eps})$,
uniformly over $1\le j,k\le N$.
Suppose that
\begin{align*}
f(\zvec) &= \exp\biggl(
   -\mwA N\sum_{j=1}^N z_j^2
   + \sum_{j=1}^{N} \mwa_j z_j^2 + N \sum_{j=1}^{N} \mwB_j z_j^3
   + \sum_{j,k=1}^N \mwC_{jk}z_j z_k^2
\\
&\kern14mm{}
   + N \sum_{j=1}^N \mwE_j z_j^4
   + \sum_{j,k=1}^N \mwF_{jk} z_j^2 z_k^2
  + \sum_{j=1}^{N} \mwJ_j z_j + \delta(\zvec) \biggr)
\end{align*}
is integrable for $\zvec=(z_1,z_2,\ldots,z_N)\in U_N$
and\/ $\delta(N)=\max_{\zvec\in U_N} \abs{\delta(\zvec)} = o(1)$,
where
\[
U_N = \bigl\{ \zvec \subseteq\Reals^N\bigm| \abs{z_j} \le N^{-1/2+\hat\eps}
  \text{ for\/ $1\le j\le N$}\bigr\},
\]
where $\hat\eps=\hat\eps(N)$ satisfies
$\eps''\le2\hat\eps\le\eps'''$.
Then, provided the $O(\,)$ term in the following
converges to zero,
\[
 \int_{U_N}f(\zvec)\,d\zvec
  = \biggl(\frac{\pi}{\mwA N}\biggr)^{\!N/2}
   \!\exp\(
     \Thetait_1+\Thetait_2 + O\((N^{-\Deltait}+\delta(N)) \mwZ\)
        \),
\]
where
\begin{align*}
  \Thetait_1 &= \frac{1}{2\mwA N}\sum_{j=1}^{N}\mwa_j
     + \frac{1}{4\mwA^2N^2}\sum_{j=1}^{N}\mwa_j^2
     + \frac{15}{16\mwA^3N}\sum_{j=1}^{N}\mwB_j^2
     + \frac{3}{8\mwA^3N^2}\sum_{j,k=1}^{N}\mwB_j\mwC_{jk} \\
  &\quad{} + \frac{1}{16\mwA^3N^3}\sum_{j,k,\ell=1}^N \mwC_{jk}\mwC_{j\ell}
    + \frac{3}{4\mwA^2N}\sum_{j=1}^N E_j
    + \frac{1}{4\mwA^2N^2}\sum_{j,k=1}^N F_{jk}\,,
  \displaybreak[0]\\
  \Thetait_2 &= \frac{1}{6\mwA^3N^3}\sum_{j=1}^N \mwa_j^3
 + \frac{3}{2\mwA^3N^2}\sum_{j=1}^N \mwa_j\mwE_j
 + \frac{45}{16\mwA^4N^2}\sum_{j=1}^N \mwa_j\mwB_j^2 \\
 &\quad{} + \frac{1}{4\mwA^3N^3}\sum_{j,k=1}^N (\mwa_j+\mwa_k) \mwF_{jk}
       + \frac{3}{4\mwA^2N}\sum_{j=1}^N \mwB_j\mwJ_j
       + \frac{1}{4\mwA^2N^2}\sum_{j,k=1}^N \mwC_{jk}\mwJ_j \\
 &\quad{} + \frac{1}{16\mwA^4N^4}\sum_{j,k,\ell=1}^N
 			(\mwa_j+2\mwa_k ) \mwC_{jk}\mwC_{j\ell}
       + \frac{3}{8\mwA^4N^3}\sum_{j,k=1}^N
                           (2\mwa_j+\mwa_k ) \mwB_j\mwC_{jk}\,,
\displaybreak[0]\\
  \mwZ &= \exp\biggl(
       \frac{1}{4\mwA^2N^2}\sum_{j=1}^{N}\Im(\mwa_j)^2
     + \frac{15}{16\mwA^3N}\sum_{j=1}^{N}\Im(\mwB_j)^2 \\
     &\kern15mm{}
     + \frac{3}{8\mwA^3N^2}\sum_{j,k=1}^{N}\Im(\mwB_j)\Im(\mwC_{jk})
     + \frac{1}{16\mwA^3N^3}\sum_{j,k,\ell=1}^N \Im(\mwC_{jk})\Im(\mwC_{j\ell})
    \biggr).
\end{align*}
\end{thm}

\begin{proof}
Our method of proof will be integration over one variable at a time.
This method is conceptually simple but technically challenging.
Assistance from a computer-algebra system is recommended.

Let $H_{j_1,j_2,\ldots,j_k}$ be a functions of $N$ for each
$1\le j_1,j_2,\ldots,j_k \le N$
and let $p_1,p_2,\ldots,p_k$ be non-negative integers.
Let $1\le j\le N+1$.  Define the generalized moment
\[\Mom_j\(H_{j_1,\ldots,j_k} |\, p_1,\ldots,p_k\)(\zvec)
   = \sum_{j_1,\ldots, j_k}
  H_{j_1,\ldots,j_k}\,z_{j_1}^{p_1}\cdots z_{j_k}^{p_k},\]
where the summation is over
\[\bigl\{ (j_1,\ldots,j_k) \bigm|
   \card{\{j_1,\ldots,j_k\}} = k, 1\le j_i\le j-1\text{ if }p_i=0,
   j\le j_i\le N\text{ if }p_i > 0 \bigr\}.\]
We will customarily omit the argument $\zvec$ as it will
be clear from the context.
Note that the indices $j_1,\ldots,j_k$ are
reserved to this notation and always index the position their name suggests;
for example
\[\Mom_j(\alpha_{j_2}|\,3,0)
 = \sum_{\substack{j\le j_1\le N\\[0.1ex]1\le j_2\le j-1}}
  \alpha_{j_2} z_{j_1}^3.\]
We also need the defective moment
$\Mom'_j(H_{j_1,\ldots,j_k} |\, p_1,p_2,\ldots,p_k)$
which is the same as\\ $\Mom_j(H_{j_1,\ldots,j_k} |\, p_1,p_2,\ldots,p_k)$
except that the index value $j$ is forbidden; that is, the condition
$j_1,\ldots,j_k\ne j$ is added to the domain of summation.

Some properties of these moments that we require are listed below.
Assume that $\zvec\in U_N$.  Then
\begin{align}
  \bigl| \Mom_j(H_{j_1,\ldots,j_k} |\, p_1,\ldots,p_k) \bigr|
     &\le  \max \,\abs{H_{j_1,\ldots,j_k}} \;
        N^{k+(-1/2+\bar\eps)(p_1+\cdots+p_k)},
\label{Mom1} \\[1ex]
 \Mom_1(H_{j_1,\ldots,j_k} |\, p_1,\ldots,p_k) &= 0
    \text{ if $p_i=0$ for any $i$} \label{Mom2},\\[0.5ex]
 \Mom_{N+1}(H_{j_1,\ldots,j_k} |\, p_1,\ldots,p_k) &= 0
    \text{ if $p_i>0$ for any $i$}, \label{Mom3}
       \displaybreak[0]\\[1ex]
 \Mom_j(H_{j_1,\ldots,j_k} |\, p_1,\ldots,p_k)
     &= \Mom'_j(H_{j_1,\ldots,j_k} |\, p_1,\ldots,p_k)\notag\\
   &~{}+ \sum_{i \,|\, p_i>0} z_j^{p_i} \,
     \Mom'_j(H_{j_1,\ldots,j_{i-1},j,j_i,\ldots,j_{k-1}}
         |\, p_1,\ldots,p_{i-1},p_{i+1},\ldots,p_k),
        \label{Mom4}\\[0.5ex]
 \Mom_{j+1}(H_{j_1,\ldots,j_k} |\, p_1,\ldots,p_k)
     &= \Mom'_j(H_{j_1,\ldots,j_k} |\, p_1,\ldots,p_k)\notag\\
   &~{}- \sum_{i \,|\, p_i=0}
     \Mom'_j(H_{j_1,\ldots,j_{i-1},j,j_i,\ldots,j_{k-1}}
         |\, p_1,\ldots,p_{i-1},p_{i+1},\ldots,p_k).
       \label{Mom5}
\end{align}
The last two equalities require $j\le N$.

The product of generalized moments
$\Mom_j(P_{j_1,\ldots,j_k} |\, p_1,\ldots,p_k)$
and $\Mom_j(Q_{j_1,\ldots,j_\ell} |\, q_1,\ldots,q_\ell)$
can be written as a sum of generalized moments.
Define $\Phi$ to be the set of injections
$\phi : \{1,2,\ldots,\ell\} \to \{1,2,\ldots,k+\ell\}$ such that
(a)~$\phi\( \{1,2,\ldots,\ell\} \) \cup \{1,2,\ldots,k\}
  = \{ 1,2,\ldots,\abs\phi\}$
for some integer $\abs\phi$ depending on~$\phi$,
(b)~for $1\le i<j\le\ell$, if $\phi(i),\phi(j)>k$
then $\phi(i)<\phi(j)$, and
(c) for $1\le i\le\ell$, $q_i=0
 \Leftrightarrow (\phi(i)>k\text{ or } p_{\phi(i)}=0)$.
For $\phi\in\Phi$ and $1\le i\le\abs\phi$, define $r_i=p_i+q_{\phi^{-1}(i)}$,
where the first term is omitted if $i>k$ and the second term is
omitted if $i$ is not in the range of~$\phi$. Then
\begin{align}
\begin{split}\label{Mom6}
   \Mom_j(P_{j_1,\ldots,j_k} |\, p_1,\ldots,p_k) &\;
     \Mom_j(Q_{j_1,\ldots,j_\ell} |\, q_1,\ldots,q_\ell)\\
   &\kern15mm{}= \sum_{\phi\in\Phi}
      \Mom_j(P_{j_1,\ldots,j_k} Q_{j_{\phi(1)},\ldots,j_{\phi(\ell)}}
       |\, r_1,\ldots,r_{\abs\phi}).
\end{split}
\end{align}
For example,
\[\Mom_j(\alpha_{j_1,j_2}|\,0,2)\,\Mom_j(\beta_{j_1}|\,3) =
   \Mom_j(\alpha_{j_1,j_2}\beta_{j_3}|\,0,2,3) +
   \Mom_j(\alpha_{j_1,j_2}\beta_{j_2}|\,0,5),\]
where the two terms correspond to the injections $\phi(1)=3$ and $\phi(1)=2$.
Exactly the same formula holds for defective moments.

For $1\le j\le N+1$, define
\def\gco#1#2{\Gammait_{#1}}
\begin{align*}
 F_j(\zvec) &= \Momj(-\mwA N+\mwa_{j_1}|\,2)
+ \Momj(\mwB_{j_1}N+\mwC_{j_1j_1}|\,3)
+ \Momj(\mwE_{j_1} N+\mwF_{j_1j_1}|\,4)
+ \Momj(\mwJ_{j_1}|\,1)
\\[0.2ex]&\quad{}
+ \Momj(\mwC_{j_2j_1}|\,2,1)
+ \Momj(\mwF_{j_2j_1}|\,2,2)
+ \Momj(\gco{0}{j_1}|\,0)
+ \Momj(\gco{0,0}{j_1j_2}|\,0,0)
+ \Momj(\gco{1,0}{j_1j_2}|\,1,0)
\\[0.4ex]&\quad{}
+ \Momj(\gco{2,0}{j_1j_2}|\,2,0)
+ \Momj(\gco{0,0,0}{j_1j_2j_3}|\,0,0,0)
+ \Momj(\gco{1,0,0}{j_1j_2j_3}|\,1,0,0)
+ \Momj(\gco{1,1,0}{j_1j_2j_3}|\,1,1,0)
\\[0.4ex]&\quad{}
+ \Momj(\gco{2,0,0}{j_1j_2j_3}|\,2,0,0)
+ \Momj(\gco{2,1,0}{j_1j_2j_3}|\,2,1,0)
+ \Momj(\gco{2,2,0}{j_1j_2j_3}|\,2,2,0)
+ \Momj(\gco{1,0,0,0}{j_1j_2j_3j_4}|\,1,0,0,0)
\\[0.4ex]&\quad{}
+ \Momj(\gco{1,1,1,0}{j_1j_2j_3j_4}|\,1,1,1,0)
+ \Momj(\gco{2,1,0,0}{j_1j_2j_3j_4}|\,2,1,0,0)
+ \Momj(\gco{2,2,1,0}{j_1j_2j_3j_4}|\,2,2,1,0),
\end{align*}
where
\begin{align*}
\gco{0}{j_1} &=
\frac{\mwa_{j_1}}{2\mwA N}
+\frac{\mwa_{j_1}^2}{4\mwA^2 N^2}
+\frac{\mwa_{j_1}^3}{6\mwA^3 N^3}
+\frac{3\mwJ_{j_1}\mwB_{j_1}}{4\mwA^2 N}
+\frac{3\mwE_{j_1}}{4\mwA^2 N}
+\frac{15\mwB_{j_1}^2}{16\mwA^3 N}
\\&\quad{}
+\frac{45\mwa_{j_1}\mwB_{j_1}^2}{16\mwA^4 N^2}
+\frac{3\mwa_{j_1}\mwE_{j_1}}{2\mwA^3 N^2},
\\
\gco{0,0}{j_1j_2} &=
\frac{\mwF_{j_2j_1}}{4\mwA^2 N^2}
+\frac{\mwC_{j_2j_1}\mwJ_{j_2}}{4\mwA^2 N^2}
+\frac{3\mwC_{j_2j_1}\mwB_{j_2}}{8\mwA^3 N^2}
+\frac{(\mwa_{j_1}+\mwa_{j_2})\mwF_{j_2j_1}}{4\mwA^3 N^3}
+\frac{3(\mwa_{j_1}+2\mwa_{j_2})\mwC_{j_2j_1}\mwB_{j_2}}{8\mwA^4 N^3},
\displaybreak[0]\\
\gco{1,0}{j_1j_2} &=
\frac{\mwC_{j_1j_2}}{2\mwA N}
+\frac{\mwa_{j_2}\mwC_{j_1j_2}}{2\mwA^2 N^2}
+\frac{45\mwC_{j_1j_2}\mwB_{j_2}^2}{16\mwA^4 N^2}
+\frac{\mwa_{j_2}^2\mwC_{j_1j_2}}{2\mwA^3 N^3}
+\frac{3\mwC_{j_1j_2}\mwE_{j_2}}{2\mwA^3 N^2},
\displaybreak[0]\\
\gco{2,0}{j_1j_2} &=
\frac{3\mwC_{j_2j_1}\mwB_{j_2}}{4\mwA^2 N}
+\frac{\mwF_{j_1j_2}+\mwF_{j_2j_1}}{2\mwA N}
+\frac{\mwC_{j_2j_1}\mwJ_{j_2}}{2\mwA N}
+\frac{\mwa_{j_2}(\mwF_{j_1j_2}+\mwF_{j_2j_1})}{2\mwA^2 N^2}
+\frac{3\mwa_{j_2}\mwC_{j_2j_1}\mwB_{j_2}}{2\mwA^3 N^2},
\displaybreak[0]\\
\gco{0,0,0}{j_1j_2j_3} &=
\frac{\mwC_{j_3j_1}\mwC_{j_3j_2}}{16\mwA^3 N^3}
+\frac{(2\mwa_{j_2}+\mwa_{j_3})\mwC_{j_3j_1}\mwC_{j_3j_2}}{16\mwA^4 N^4},
\displaybreak[0]\\
\gco{1,0,0}{j_1j_2j_3} &=
\frac{(\mwF_{j_2j_3}+\mwF_{j_3j_2})\mwC_{j_1j_2}}{4\mwA^3 N^3}
+\frac{3\mwC_{j_1j_3}\mwC_{j_3j_2}\mwB_{j_3}}{4\mwA^4 N^3}
+\frac{3\mwC_{j_3j_2}\mwC_{j_1j_2}\mwB_{j_3}}{8\mwA^4 N^3},
\displaybreak[0]\\
\gco{1,1,0}{j_1j_2j_3} &=
\frac{\mwC_{j_1j_3}\mwC_{j_2j_3}}{4\mwA^2 N^2}
+\frac{\mwa_{j_3}\mwC_{j_1j_3}\mwC_{j_2j_3}}{2\mwA^3 N^3},
\displaybreak[0]\\
\gco{2,0,0}{j_1j_2j_3} &=
\frac{(\mwa_{j_2}+\mwa_{j_3})\mwC_{j_3j_2}\mwC_{j_3j_1}}{4\mwA^3 N^3}
+\frac{\mwC_{j_3j_1}\mwC_{j_3j_2}}{4\mwA^2 N^2},
\displaybreak[0]\\
\gco{2,1,0}{j_1j_2j_3} &=
\frac{3\mwC_{j_2j_3}\mwC_{j_3j_1}\mwB_{j_3}}{2\mwA^3 N^2}
+\frac{(\mwF_{j_1j_3}+\mwF_{j_3j_1})\mwC_{j_2j_3}}{2\mwA^2 N^2},
\displaybreak[0]\\
\gco{2,2,0}{j_1j_2j_3} &=
\frac{\mwC_{j_3j_1}\mwC_{j_3j_2}}{4\mwA N}
+\frac{\mwa_{j_3}\mwC_{j_3j_1}\mwC_{j_3j_2}}{4\mwA^2 N^2},
\displaybreak[0]\\
\gco{1,0,0,0}{j_1j_2j_3j_4} &=
\frac{\mwC_{j_1j_4}\mwC_{j_4j_2}\mwC_{j_4j_3}}{16\mwA^4 N^4}
+\frac{\mwC_{j_4j_3}\mwC_{j_4j_2}\mwC_{j_1j_2}}{8\mwA^4 N^4},
\displaybreak[0]\\
\gco{1,1,1,0}{j_1j_2j_3j_4} &=
\frac{\mwC_{j_1j_4}\mwC_{j_2j_4}\mwC_{j_3j_4}}{6\mwA^3 N^3},
\\
\gco{2,1,0,0}{j_1j_2j_3j_4} &=
\frac{\mwC_{j_2j_4}\mwC_{j_4j_1}\mwC_{j_4j_3}}{4\mwA^3 N^3}
+\frac{\mwC_{j_4j_1}\mwC_{j_4j_3}\mwC_{j_2j_3}}{4\mwA^3 N^3},
\\
\gco{2,2,1,0}{j_1j_2j_3j_4} &=
\frac{\mwC_{j_3j_4}\mwC_{j_4j_1}\mwC_{j_4j_2}}{4\mwA^2 N^2}\,.
\end{align*}

Note that $F_j(\zvec)$ is independent of $z_i$ for $i < j$.
The key properties of $F_j(\zvec)$ for $\zvec\in U_n$ are
\begin{align}
   f(\zvec) &= \exp\( F_1(\zvec) + \delta(\zvec) \), \label{Mom7} \\[1ex]
   \int_{-N^{-1/2+\hat\eps}}^{N^{-1/2+\hat\eps}} \exp\(F_j(\zvec)\)\,dz_j
     &= \sqrt{\!\frac{\pi}{\mwA N}} \exp\(F_{j+1}(\zvec) + O(N^{-1-\Deltait})\)
        \qquad(j \le N). \label{Mom8}
\end{align}

Equation \eqref{Mom7} is easily seen after applying \eqref{Mom2} to eliminate
most of the terms.  Proof of \eqref{Mom8} requires a tedious calculation
which we now outline.

First, apply \eqref{Mom4} to make explicit the dependence of $F_j(\zvec)$
on~$z_j$ (a polynomial of degree~4).  Then expand
\begin{equation}\label{Fjj}
  \exp\(F_j(\zvec)\) = \exp\(R_0(\zvec)\)
      \exp\(-\mwA N z_j^2\)
       \( 1 + R_1(\zvec)z_j + R_2(\zvec)z_j^2 + \cdots + O(N^{-1-\Deltait})\),
\end{equation}
where each $R_i(\zvec)$ is independent of $z_j$ and
contains defective moments only.
As seen by applying~\eqref{Mom1},
only a finite number of terms are required to achieve the requested
error term.  The factor $\exp(-\mwA N z_j^2)$ comes from the
first term $\Momj(-\mwA N|\,2)$ of $F_j(\zvec)$.  Products of moments
that occur need to be rewritten as sums using~\eqref{Mom6}.

Next, integrate \eqref{Fjj} over $z_j$ using
\[\int_{-N^{-1/2+\hat\eps}}^{N^{-1/2+\hat\eps}} z^{2k}e^{-\mwA Nz^2}\,dz
   = \frac{(2k)!}{k!\,(4\mwA N)^k} \sqrt{\!\frac{\pi}{\mwA N}}
       \(1 + O(\exp(-c N^{2\hat\eps-\eps'}))\),\]
for fixed $k\ge 0$, for some $c>0$.  Here we have used the assumptions
that $\mwA=\Omega(N^{-\eps'})$ and $\eps'<\eps''<2\hat\eps$.
The result of the integration has the form
\[ \sqrt{\!\frac{\pi}{\mwA N}}
     \exp\(R_0(\zvec)\) \( 1 + S(\zvec) \)
     = \sqrt{\!\frac{\pi}{\mwA N}} \exp\(R_0(\zvec) + \log(1+S(\zvec))\).\]
Since $S(\zvec)=o(1)$
(in fact $S(\zvec)=O(N^{-1/2+k(\eps''+\bar\eps)})$
for some~$k$), we can expand the logarithm using~\eqref{Mom1}
again to limit the expansion to finitely many terms.  Finally,
apply~\eqref{Mom5} to rewrite the defective moments in terms of
ordinary generalized moments.
The result is the right side of~\eqref{Mom8}.

\smallskip

If all the coefficients $\mwa_j, \mwB_j, \mwC_{jk}, \mwE_j, \mwF_{jk},
\mwJ_j$ were real,
we could apply~\eqref{Mom7} and~\eqref{Mom8} immediately to find that
\begin{equation}\label{Mom9}
\int_{U_N} f(\zvec)\,d\zvec =
\biggl( \frac{\pi}{\mwA N} \biggr)^{\!N/2}
   \exp\( F_{N+1} + O(\delta(N)+N^{-\Deltait})\),
\end{equation}
noting that $F_{N+1}(\zvec)$ is independent of~$\zvec$.

When the coefficients are complex, we must take more care.
Equation~\eqref{Mom7} only allows us to write
\begin{equation}\label{Mom10}
   \int_{U_N} f(\zvec)\,d\zvec
     = \int_{U_N} \exp\(F_1(\zvec)\) \,d\zvec
        + O(\delta(N)) \int_{U_N} \bigl|\exp\(F_1(\zvec)\)\bigr| \,d\zvec.\\
\end{equation}

Let $F_j^*(\zvec)$ be the same as $F_j(\zvec)$ except that the
coefficients  $\mwa_j, \mwB_j, \mwC_{jk}, \mwE_j, \mwF_{jk},
\mwJ_j$ are all replaced
by their real parts.  Clearly
\[\bigl| \exp(F_1(\zvec)) \bigr| = \exp(F_1^*(\zvec)),\]
and so, as in \eqref{Mom9},
\begin{equation}\label{Mom11}
\int_{U_N} \bigl| \exp(F_1(\zvec)) \bigr| \,d\zvec =
  O(1) \biggl( \frac{\pi}{\mwA N} \biggr)^{\!N/2}\exp(F_{N+1}^*).
\end{equation}
{}From \eqref{Mom8} we have for $1\le j\le N$ that
\begin{align*}
\bigl| \exp(F_{j+1}(\zvec)) \bigr|
   &= \biggl(\frac{\pi}{\mwA N}\biggr)^{\!-1/2}\(1+O(N^{-1-\Deltait})\)
              \,\biggl| \int \exp(F_j(\zvec))\,dz_j \Biggr| \\
   &\le \biggl(\frac{\pi}{\mwA N}\biggr)^{\!-1/2}\(1+O(N^{-1-\Deltait})\)
              \int \,\bigl|\exp(F_j(\zvec)) \bigr|\,dz_j,
\end{align*}
so we have by induction starting with \eqref{Mom11} that
\[
  \int \,\bigl| \exp(F_{j+1}(\zvec)) \bigr|\,dz_{j+1}\cdots dz_N
     \le O(1) \biggl(\frac{\pi}{\mwA N}\biggr)^{\!(N-j)/2}
         \exp(F_{N+1}^*).
\]

Returning to \eqref{Mom8}, we find that, for $1\le j\le N$,
\begin{align*}
 \int \exp\(F_j(\zvec)\)\,dz_j\cdots dz_N
   &= \sqrt{\!\frac{\pi}{\mwA N}}
       \int \exp\(F_{j+1}(\zvec)+O(N^{-1-\Deltait})\)\,
                     dz_{j+1}\cdots dz_N \\[0.4ex]
   &= \sqrt{\!\frac{\pi}{\mwA N}}
       \int \exp\(F_{j+1}(\zvec)\)\,dz_{j+1}\cdots dz_N \\
    &\qquad{}+  O(N^{-1-\Deltait}) \sqrt{\!\frac{\pi}{\mwA N}}
         \int \big|\exp(F_{j+1}(\zvec))\bigr|\,dz_{j+1}\cdots dz_N \\
      &= \sqrt{\!\frac{\pi}{\mwA N}}
       \int \exp\(F_{j+1}(\zvec)\)\,dz_{j+1}\cdots dz_N \\
    &\qquad{}+  O(N^{-1-\Deltait})
         \biggl(\frac{\pi}{\mwA N}\biggr)^{\!(N-j+1)/2}
         \exp(F_{N+1}^*).
\end{align*}
By induction on $j$, this gives
\[\int_{U_N} \exp\(F_1(\zvec)\)\,d\zvec =
   \biggl(\frac{\pi}{\mwA N}\biggr)^{\!N/2} \(\exp(F_{N+1})
     + O(N^{-\Deltait})\exp(F_{N+1}^*)\),\]
which, together with \eqref{Mom10} and \eqref{Mom11} gives
\begin{align}
\int_{U_N} f(\zvec)\,d\zvec &=
   \biggl(\frac{\pi}{\mwA N}\biggr)^{\!N/2} \(\exp(F_{N+1})
     + O(N^{-\Deltait}+\delta(N))\exp(F_{N+1}^*)\)\notag \\
  &=  \biggl(\frac{\pi}{\mwA N}\biggr)^{\!N/2}
       \exp\( F_{N+1} + O(N^{-\Deltait}+\delta(N))\mwZ \),
          \label{Mom12}
\end{align}
where $\mwZ =  \exp\(F_{N+1}^*-\Re(F_{N+1})\)$ and the
last line is valid if $\(N^{-\Deltait}+\delta(N)\)\mwZ=o(1)$.

Applying \eqref{Mom3} to the definition of $F_j(\zvec)$, we find that
\[F_{N+1} = \Thetait_1 + \Thetait_2 + O(N^{-\Deltait}),\]
from which it follows that $\mwZ$ has the value in the
theorem statement to within a multiplied constant.
Also note that $\mwZ\ge 1$, which is easiest to see by
noting that the argument of the exponential is a non-negative
quadratic form for each~$j$.
The theorem now follows from~\eqref{Mom12}.
\end{proof}

\nicebreak

\end{document}